# The Spinor Representation of Minimal Surfaces

Rob Kusner
Mathematics Department
University of Massachusetts at Amherst

Nick Schmitt
Center for Geometry, Analysis, Numerics and Graphics
University of Massachusetts at Amherst

# Contents

# Abstract

The spinor representation is developed and used to investigate minimal surfaces in $\mathbb{R}^3$ with embedded planar ends. The moduli spaces of planar-ended minimal spheres and real projective planes are determined, and new families of minimal tori and Klein bottles are given. These surfaces compactify in $S^3$ to yield surfaces critical for the Möbius invariant squared mean curvature functional $W$. On the other hand, all $W$-critical spheres and real projective planes arise this way. Thus we determine at the same time the moduli spaces of $W$-critical spheres and real projective planes via the spinor representation.

# Introduction

In this paper we investigate the interplay between spin structures on a Riemann surface $M$ and immersions of $M$ into three-space. Here, a spin structure is a complex line bundle $S$ over $M$ such that $S \otimes S$ is the holomorphic (co)tangent bundle $T(M)$ of $M$. Thus we may view a section of a $S$ as a "square root" of a holomorphic 1-form on $M$. Using this notion of spin structure, in the first part of this paper we develop the notion of the *spinor representation of a surface in space*, based on an observation of Dennis Sullivan [27]. The classical Weierstrass representation is

$$(g, \eta) \longrightarrow \text{Re} \int (1 - g^2, i(1 + g^2), 2g)\eta,$$

where $g$ and $\eta$ are respectively a meromorphic function and one-form on the underlying compact Riemann surface. The spinor representation (Theorem 5) is

$$(s_1, s_2) \longrightarrow \text{Re} \int (s_1^2 - s_2^2, i(s_1^2 + s_2^2), 2s_1 s_2),$$

where $s_1$ and $s_2$ are meromorphic sections of a spin structure $S$. Either representation gives a (weakly) conformal harmonic map $M \to \mathbb{R}^3$, which therefore parametrizes a (branched) minimal surface.

One feature of the spinor representation is that fundamental topological information, such as the regular homotopy class of the immersion, can be read off directly from the analytic data (Theorem 6). In fact, for the special case where the Riemann surface $M$ is hyperelliptic, we are able to give an explicit calculation of the Arf invariant for the immersion (Theorem 8); the Arf invariant distinguishes whether or not an immersion of an orientable surface is regularly homotopic to an embedding. We also consider in Part I the spinor representation for nonorientable minimal surfaces in terms of a lifting to the orientation double cover (Theorem 11). This is sufficient for constructing examples later in the paper, but is less satisfying theoretically. In a future paper, we plan to consider the general case from the perspective of "pin" structures, and also give a more direct differential geometric treatment of the Arf invariant.

The second part of this paper focuses on general properties of minimal surfaces with embedded planar ends from the viewpoint of the spinor representation. It is well-known (see [2], [13], [14]) that such surfaces conformally compactify to give extrema for the squared mean curvature integral $W = \int H^2 dA$ popularized by Willmore. Conversely, for genus zero, all $W$-critical surfaces arise this way [2].

Using the spinor representation to study these special minimal surfaces has the computational advantage of converting certain quadratic conditions to linear ones. This is carried out in Part II, where we refine the tools we need. In fact, associated to a spin structure $S$ on a closed orientable Riemann surface $M$ is a vector space $\mathcal{K}$ of sections of $S$ such that pairs of independent sections $(s_1, s_2)$ from $\mathcal{K}$ form the spinor representations of all the minimal immersions of $M$ with embedded planar ends (Theorem 13). Thus the problem of finding all these immersions is reduced to an algebraic problem (Theorem 15). In order to better understand $\mathcal{K}$, a skew-symmetric bilinear form $\Omega$ is defined from whose kernel $\mathcal{K}$ is computable (Theorem 17).

The third (and final) part of this paper is devoted to the construction of examples and to classification results. Specifically, for a given finite topological type of surface, we want to explore the moduli space $\mathcal{M}$ of immersed minimal surfaces (up to similarity) of this type with embedded planar ends: the dimension and topology of $\mathcal{M}$, convergence to degenerate cases (that is, the natural closure of $\mathcal{M}$), and examples with special symmetry (which correspond to singular points of $\mathcal{M}$). The tools mentioned above permit the broad outline of a solution, but require ingenuity to apply in particular cases. For example, the form $\Omega$ allows the moduli space to be expressed as a determinantal variety which determines how the location of the ends can vary along the Riemann surface $M$. However, this determinantal variety is only computable when the number of ends is small. Furthermore, the basic tools, being algebraic geometrical, ignore the real analytic problems of removing periods and branch points. The latter require much subtler and often *ad hoc* methods.

Previously known results concerning genus zero minimal surfaces with embedded planar ends include the following:
- examples have been found for 4, 6, and every $n \geq 8$ ends [2], [14], [23];
- there are no immersed examples with 3, 5, and 7 ends [3];
- the moduli spaces for immersed spheres with 4 and 6 ends, and projective planes with 3 ends have been determined [3].

In Part III our new theorems include the following:
- a new proof of the non-existence of examples with 3, 5 and 7 ends is given using the skew-symmetric form $\Omega$ (Theorem 18);
- the moduli space for $2p$ ends ($2 \leq p \leq 7$) is shown to be $4(p-1)$-dimensional (Theorem 21);
- the point which compactifies the moduli space of projective planes with 3 ends is proved to be a Möbius strip, and all symmetries of these surfaces are found (Theorem 25).

A recent result concerning genus one is the construction in [5] of examples with four embedded planar ends, assuming a rectangular lattice. We give further results:
- there are no three-ended tori (Theorem 26);
- there is a real two-dimensional family of four-ended immersed examples on each conformal type of torus (Theorem 27);
- there exists an immersed Klein bottle with four ends (Theorem 29).

For higher genus, the general methods we have developed here also yield (possibly branched) minimal immersions with embedded planar ends, but it becomes more and more difficult to determine precisely when branch points are absent or periods vanish: we again postpone this case to a future paper.

Most of the theorems presented here were worked out while we visited the Institute for Advanced Study during the 1992 Fall term, and were first recorded in [26]. It is a pleasure to thank the School of Mathematics at the Institute for its hospitality, as well as Sasha Bobenko, Peter Norman and Dennis Sullivan for their comments and interest. In particular, we should mention that Bobenko has recently announced some related results for constant mean curvature surfaces (*Surfaces in terms of 2 by 2 matrices: Old and new integrable cases,* in: A. Fordy and J. Wood, *Harmonic Maps and Integrable Systems*: Vieweg, 1994).

# Part I
**Spinors, Regular Homotopy Classes and the Arf Invariant**

## 1 The spinor representation

The notion of a spin structure is developed and used to describe the spinor representation of a surface in space. Section 3 defines a "quadratic form" which can be used to completely classify the spin structures on a surface, and section 4 computes coordinates for the unique spin structure on the Riemann sphere. In the next two sections, the spinor representation of a surface is explained and related to the regular homotopy class of the surface. Section 7 shows equivalent characterizations of spin structures, the most useful of which will be that of representing spin structures by holomorphic differentials. These differentials are computed on hyperelliptic Riemann surfaces. Section 9 takes up the question of group action on spinors, and computes the group which performs Euclidean similarity transformations. Two surfaces which are transforms of each other under the action of this group are considered to be the same. The final two sections discuss briefly the technicalities of periods and nonorientable surfaces.

## 2 Spin structures and spin manifolds

A spin structure on an $n$-dimensional (spin) manifold $M$ is a certain two-sheeted covering map of the SO($n$)-frame-bundle on $M$ to a Spin($n$)-bundle (see [20], [17]). When $n = 2$, this notion of spin structure may easily be reduced to the following definition in terms of a quadratic map between complex line bundles, as the figure below depicts:

$$\begin{array}{ccc} S & \xrightarrow{\mu} & T(M) \\ & \searrow & \downarrow \\ & & M \end{array}$$

Figure 1: Spin structure

<u>DEFINITION 1.</u> *A* spin structure *on a Riemann surface $M$ is a complex line bundle $S$ over $M$ together with a smooth surjective fiber-preserving map $\mu : S \longrightarrow T(M)$ to* *the holomorphic (co)tangent bundle $T(M)$ satisfying*

(2.1) $$\mu(\lambda s) = \lambda^2 \mu(s)$$

*for any section $s$ of $S$.*

Two spin structures $(S, \mu)$ and $(S', \mu')$ on a Riemann surface $M$ are *isomorphic* if there is a line bundle isomorphism $\delta : S \longrightarrow S'$ for which $\mu = \mu'\delta$. Hence two spin structures may be isomorphic as line bundles and yet not be isomorphic as spin structures. The number of nonisomorphic spin structures on a Riemann surface $M$ is equal to the cardinality of $H^1(M, \mathbb{Z}_2)$. (This count remains true for spin manifolds in general: see [20].) In particular, if $M$ is a closed Riemann surface of genus $g$, there are $2^{2g} = \#H^1(M, \mathbb{Z}_2)$ such structures on $M$.

On an annulus $A = \{r_1 < z < r_2\}$ there are exactly two nonisomorphic spin structures, which can be given explicitly as follows. The tangent bundle $T(A)$ may be identified with $A \times \mathbb{C}$ by means of the global trivialization

$$a \left.\frac{\partial}{\partial z}\right|_p \mapsto (p, a).$$

Let $S_0 = S_1 = A \times \mathbb{C}$ and define maps $\mu_k : S_k \longrightarrow T(A)$ for $k = 0, 1$ by

$$\mu_0(z, w) = (z, w^2),$$
$$\mu_1(z, w) = (z, zw^2).$$

Then $(S_k, \mu_k)$ are spin structures on $A$ since $\mu_k$ satisfies the condition (2.1). Though $S_0$ and $S_1$ are isomorphic line bundles over $A$, they are nonisomorphic spin structures. For if $S_0$ and $S_1$ were isomorphic spin structures with bundle isomorphism $\delta : S_0 \longrightarrow S_1$ satisfying $\mu_0 = \mu_1\delta$, then $\delta$ would be of the form $(z, w) \mapsto (z, f(z, w))$. Then $w^2 = zf^2$, implying that $z$ has a consistent square root on $\mathbb{C}^*$, which is impossible.

## 3  The quadratic form associated to a spin structure

In this section, the Riemann surface $M$, its holomorphic (co)tangent bundle, and the spin structure are replaced with the corresponding real manifold and real vector bundles. In particular, all vector fields in this section are *real* vector fields.

To each spin structure $S$ on the Riemann surface $M$ we associate a $\mathbb{Z}_2$-valued quadratic form

$$q : H_1(M, \mathbb{Z}_2) \longrightarrow \mathbb{Z}_2.$$

To say that $q$ is quadratic means that for all $\alpha_1, \alpha_2 \in H_1(M, \mathbb{Z}_2)$ we have

$$q(c_1 + c_2) = q(c_1) + q(c_2) + c_1 \cdot c_2.$$

where $c_1 \cdot c_2$ denotes the mod 2 intersection number of $c_1$ with $c_2$.

To define $q(c)$, let $\alpha : S^1 \longrightarrow M$ be a smooth embedded representative of $c$ (the existence of such an $\alpha$ follows from results in [19]). Let $v$ be a smooth vector field along $\alpha$ which lifts to a section of the spin structure along $\alpha$, and let $w(\alpha, v)$ denote the total turning number, mod 2, of the derivative vector $\alpha'$ against $v$ along $\alpha$. Define $q(c) = w_v(c) + 1$.

To show that $q$ is quadratic, the following technical lemma is stated without proof. (A *Jordan trail* is a closed tracing along a curve which tracing does not cross itself. The existence of the Jordan trail is assured in [12].)

LEMMA 2.
(i) Let $\alpha : S^1 \longrightarrow M$ be an immersion, and let $a$ be the number of self-crossing points of $\alpha$. Let $v$ be a smooth non-zero vector field along $\alpha(S^1)$ on $M$. Let $\beta$ be a Jordan trail for $\alpha$. Then

$$w_v(\beta) = w_v(\alpha) + a.$$

(ii) Let $\alpha_1$ and $\alpha_2 : S^1 \longrightarrow M$ be immersions, with a common base point $\alpha_1(t) = \alpha_2(t)$. Let $\alpha_1 * \alpha_2 : S^1 \longrightarrow M$ denote the closed curve consisting of $\alpha_1$ followed by $\alpha_2$. Let $v$ be a smooth non-zero vector field along $\alpha_1(S^1) \cup \alpha_2(S^1)$. Then

$$w_v(\alpha_1 * \alpha_2) = w_v(\alpha_1) + w_v(\alpha_2).$$

LEMMA 3. If $\alpha : S^1 \longrightarrow M$ is an embedded curve on a spin surface $M$ with spin structure $S$, and $v_1$, $v_2$ are smooth nonzero vector fields along $\alpha$, then the following are equivalent:
(i) $w(\alpha, v_1) = w(\alpha, v_2)$
(ii) $v_1$ and $v_2$ alike lift or do not lift along $\alpha$ to smooth sections of $S$.

*Proof.* We may assume $M$ is an annulus containing $\alpha(S^1)$ as the unit circle, with spin structure $S_k$ ($k = 0$ or $1$) as in section 2. Any vector field $S^1 \longrightarrow \mathbb{C}$ is of the form $t \mapsto t^p [f(t)]^2$, where $f$ is smooth and

$$p = \begin{cases} k & \text{if } v \text{ lifts,} \\ 1 - k & \text{if } v \text{ does not lift.} \end{cases}$$

Then, with $w_\alpha(h_1, h_2)$ defined as the mod 2 winding number of $h_1$ against $h_2$ (or equivalently, of $h_2/h_1$) along $\alpha$,

$$w_\alpha(v_1, v_2) = w_\alpha(t^p[f_1(t)]^2, t^q[f_2(t)]^2) = w_\alpha(t^p, t^q) \equiv p + q \pmod 2$$
$$= \begin{cases} 0 & \text{if } v_1, v_2 \text{ alike lift or do not lift,} \\ 1 & \text{otherwise.} \end{cases}$$

But $w_\alpha(v_1, v_2) = w(\alpha, v_1) + w(\alpha, v_2)$, and the result follows. □

THEOREM 4. *The form $q : H_1(M, \mathbb{Z}_2) \longrightarrow \mathbb{Z}_2$ defined above is well-defined, that is, independent of the choice of the vector field $v$ and the choice of embedded representative $\alpha$. Moreover, $q$ satisfies*

$$q(c_1 + c_2) = q(c_1) + q(c_2) + c_1 \cdot c_2.$$

*Proof.* Let $\alpha_0$, and $\alpha_1 : S^1 \longrightarrow M$ be embedded representatives of $c \in H_1(M, \mathbb{Z}_2)$. Let $v_0$, $v_1$ be smooth nonzero vector fields which lift along $\alpha_0$, $\alpha_1$ respectively to sections of the spin structure $S$. Let $\alpha_t$ ($t \in [0, 1]$) be a homotopy of $\alpha_0$ and $\alpha_1$. Extend $v_0$ to a smooth nonzero vector field in an annulus containing the image of $\alpha_t$.

Then $w(\alpha_t, v)$ is a continuous function of $t$, and an integer, hence it is constant. In particular,
$$w(\alpha_0, v_0) = w(\alpha_1, v).$$
But $v = v_0$ lifts along $\alpha_0$ to a smooth section of $S$. So $v$ must also lift along $\alpha_1$. But since $v_1$ also lifts along $\alpha_1$,
$$w(\alpha_1, v) = w(\alpha_1, v_1).$$
Thus
$$w(\alpha_0, v_0) = w(\alpha_1, v_1),$$
showing that $q$ is well-defined.

Now, to show $q$ is quadratic, let $\alpha_1$, $\alpha_2$ be embedded representatives of $c_1$, $c_2 \in H_1(M, \mathbb{Z}_2)$, and let

$$a = \# \text{ of self-crossing points of } \alpha_1 * \alpha_2 \equiv \alpha_1 \cdot \alpha_2 - 1 \pmod 2.$$

Let $\beta$ be a Jordan trail for $\alpha_1 * \alpha_2$. Then

$$w(\beta, v) = w(\alpha_1 * \alpha_2, v) + a = w(\alpha_1, v) + w(\alpha_2, v) + \alpha_1 \cdot \alpha_2 + 1$$

by Lemma 2(i). Hence

$$\begin{aligned} q(c_1 + c_2) &= w(\beta, v) + 1 = (w(\alpha_1, v) + 1) + (w(\alpha_2, v) + 1) + \alpha_1 \cdot \alpha_2 \\ &= q(c_1) + q(c_2) + c_1 \cdot c_2. \end{aligned}$$

□

A well-known result (see, for example, [24]) is that the equivalence class of the quadratic form $q : H_1(M, \mathbb{Z}_2) \longrightarrow \mathbb{Z}_2$ under linear changes of bases of $H_1(M, \mathbb{Z}_2)$ is determined by its *Arf invariant*

$$(3.2) \qquad \operatorname{Arf} q = \frac{1}{\sqrt{\#H}} \sum_{\alpha \in H} (-1)^{q(\alpha)},$$

where $H = H_1(M, \mathbb{Z}_2)$. The quadraticity of $q$ insures that this invariant has values in $\{+1, -1\}$. For a compact surface of genus $g$, there are $2^{2g-1} + 2^{g-1}$ spin structures for which the Arf invariant of the corresponding quadratic form is $+1$, and $2^{2g-1} - 2^{g-1}$ spin structures for which it is $-1$ (see section 8).

## 4   The spin structure on the Riemann sphere

The following description of the unique spin structure on $S^2$, as well as the spinor representation of a surface in the next section, are adapted from [27]. Identify

$$S^2 \cong [Q] = \{[z_1, z_2, z_3] \in \mathbb{CP}^2 \mid z_1^2 + z_2^2 + z_3^2 = 0\},$$

where $Q$ is the null quadric

$$Q = \{(z_1, z_2, z_3) \in \mathbb{C}^3 \mid z_1^2 + z_2^2 + z_3^2 = 0\}.$$

Then $T(S^2)$ may be identified with the restriction to $[Q]$ of the tautological line bundle

$$\operatorname{Taut}(\mathbb{CP}^2) = \{(\Lambda, x) \in \mathbb{CP}^2 \times \mathbb{C}^3 \mid x \in \Lambda\}$$

(here, $\mathbb{CP}^2$ is thought of as the lines in $\mathbb{C}^3$), so

$$(4.3) \qquad T(S^2) \cong \operatorname{Taut}(\mathbb{CP}^2)|_{[Q]} = \{(\Lambda, x) \in [Q] \times Q \mid x = 0 \text{ or } \pi(x) \in \Lambda\},$$

where $\pi : Q \longrightarrow [Q]$ is the canonical projection. Given this, the unique spin structure $\operatorname{Spin}(S^2)$ on $S^2$ may then be identified with the tautological line bundle

$$(4.4) \qquad \operatorname{Spin}(S^2) \cong \operatorname{Taut}(\mathbb{CP}^1) \cong \{(\Lambda, x) \in \mathbb{CP}^1 \times \mathbb{C}^2 \mid x \in \Lambda\},$$

with the associated mapping $\mu$ given by

$$\mu([z_1, z_2], (s_1, s_2)) = ([\sigma(z_1, z_2)], \sigma(s_1, s_2)),$$

where $\sigma : \mathbb{C}^2 \longrightarrow Q$ is the map defined by

(4.5) $$\sigma(z_1, z_2) = (z_1^2 - z_2^2, i(z_1^2 + z_2^2), 2z_1 z_2).$$

As may be checked, the map $\mu$ satisfies the conditions of Definition 1.

When $T(S^2)$ and $\mathrm{Spin}(S^2)$ are restricted respectively to their nonzero vectors and nonzero spin-vectors, they have single coordinate charts

$$\begin{aligned}\{\text{nonzero vectors in } T(S^2)\} &\longrightarrow Q \setminus \{0\} \\ \{\text{nonzero spin-vectors in } \mathrm{Spin}(S^2)\} &\longrightarrow \mathbb{C}^2 \setminus \{0\}\end{aligned}$$

defined by taking the second component in each of (4.3) and (4.4) respectively. In this case, $\mu$ may be thought of as the two-to-one covering map $\sigma : \mathbb{C}^2 \setminus \{0\} \longrightarrow Q \setminus \{0\}$.

## 5  The spinor representation of a surface

To describe the spinor representation, let $M$ be a Riemann surface with a local complex coordinate $z$, and $X : M \longrightarrow \mathbb{R}^3$ a conformal (but not necessarily minimal) immersion of $M$ into space. Since $X$ is conformal, its $z$-derivative $\partial X = \omega$ can be viewed as a null vector in $\mathbb{C}^3$, or via (4.3), as a map into the (co)tangent bundle $T(S^2)$; so with the (not necessarily meromorphic) Gauss map $g$ associated to $X$, we get the bundle map $(\omega, g)$ as in Figure 2.

$$\begin{array}{ccc} T(M) & \xrightarrow{\omega} & T(S^2) \\ \downarrow & & \downarrow \\ M & \xrightarrow{g} & S^2 \end{array}$$

Figure 2: Bundle map

The Weierstrass representation is determined by $(g, \eta)$ where $\eta$ is the (not necessarily meromorphic) differential form defined by

$$\omega = (1 - g^2, i(1 + g^2), 2g)\,\eta.$$

Conversely, given a bundle map $(\omega, g)$ of $T(M)$ into $T(S^2)$, if $\omega$ satisfies the integrability condition $\mathrm{Re}\,d\omega = 0$, then

$$X = \mathrm{Re}\int \omega : M \longrightarrow \mathbb{R}^3$$

is a (possibly periodic) immersion with Gauss map $g$.

The spinor representation of the immersion, shown in Figure 3, is obtained by lifting $\omega$ to the spin structures on $M$ and $S^2$.

$$
\begin{array}{ccc}
S & \xrightarrow{\hat{\omega}} & \mathrm{Spin}(S^2) \\
\mu \downarrow & & \downarrow \sigma \\
T(M) & \xrightarrow{\omega} & T(S^2) \\
\downarrow & & \downarrow \\
M & \xrightarrow{g} & S^2
\end{array}
$$

Figure 3: Spinor representation of a surface

THEOREM 5. *Let $M$ be a connected surface, and $(\omega, g)$ a bundle map of $T(M)$ into $T(S^2)$. Then*
  (i) *there is a unique spin structure $S$ on $M$ such that $\omega$ lifts to a bundle map $\hat{\omega} : S \longrightarrow \mathrm{Spin}(S^2)$;*
  (ii) *there are exactly two such lifts $\hat{\omega}$, and these differ only by sign.*

*Proof of (i):* Considering $\mathrm{Spin}(S^2)$ as a $\mathbb{Z}_2$-bundle on $T(S^2)$ when restricted to nonzero spin-vectors and vectors respectively, let $S$ be the (unique) pullback bundle of $\mathrm{Spin}(S^2)$ under $\omega$, and $\mu$, $\hat{\omega}$ as shown. Extend $S$, $\hat{\omega}$, and $\mu$ to include the zero spin-vectors.

*Proof of (ii):* If $\iota : \mathrm{Spin}(S^2) \longrightarrow \mathrm{Spin}(S^2)$ is the order-two deck transformation for the covering $\mathrm{Spin}(S^2) \longrightarrow T(S^2)$, then $\iota \circ \hat{\omega}$ is another map which in place of $\hat{\omega}$ makes the diagram commute. Conversely, if $\zeta : S \longrightarrow \mathrm{Spin}(S^2)$ is such a map, then for $x \in S$, $\zeta(x)$ is $\hat{\omega}(x)$ or $\iota \circ \hat{\omega}(x)$ and continuity implies that $\zeta = \hat{\omega}$ or $\iota\hat{\omega}$. □

The spinor representation is determined by the pair of sections $\hat{\omega} = (s_1, s_2)$ of $S$ related to $\omega$ by the equation

$$\omega = (s_1^2 - s_2^2, i(s_1^2 + s_2^2), 2s_1 s_2).$$

Thus the Weierstrass representation and the spinor representation are related by the equations

$$\eta = s_1^2 \quad \text{and} \quad g = s_2/s_1.$$

The case of a nonorientable $M$ is dealt with in section 11 by the taking of the spin structure on the oriented two-sheeted cover of $M$.

## 6 Regular homotopy classes and spin structures

Let $X_1, X_2 : M \longrightarrow \mathbb{R}^3$ be two immersions of a surface into space. Recall the distinction between regular homotopy equivalence of the immersions $X_1$, $X_2$, and regular homotopy equivalence of the corresponding immersed surfaces — these immersed surfaces are regularly homotopic if there is a diffeomorphism $\varphi$ of $M$ such that $X_2$ is regularly homotopic to $X_1 \circ \varphi$ — so this latter equivalence relation is coarser.

THEOREM 6. *Let $X_1, X_2 : M \longrightarrow \mathbb{R}^3$ be two immersions of a surface into space, let $S_1$, $S_2$ the spin structures induced as in Theorem 5, and let $q_1$, $q_2$ be the associated quadratic forms as in Theorem 4. Then*
  (i) *$X_1$ and $X_2$ are regularly homotopic if and only if $q_1 \equiv q_2 (\mathrm{mod}\, 2)$.*
  (ii) *The surfaces $X_1(M)$ and $X_2(M)$ are regularly homotopic if and only if $\mathrm{Arf}\, q_1 = \mathrm{Arf}\, q_2$. In particular, an immersed surface is regularly homotopic to an embedding if and only if its $\mathrm{Arf}$ invariant equals $+1$.*

*Sketch of proof.* Define $\tilde{q}(\alpha)$ as the linking number (mod 2) of the boundary curves of the image of a tubular neighborhood of $\alpha$ in $\mathbb{R}^3$. Then

$$q(\alpha) = 0 \quad \Longleftrightarrow \quad \begin{array}{c} \text{the Darboux frame along } \alpha \\ \text{is nontrivial as an element of} \\ \pi_1(\mathrm{SO}(3)) \end{array} \quad \Longleftrightarrow \quad \tilde{q}(\alpha) = 0.$$

Hence $q \equiv \tilde{q} (\mathrm{mod}\, 2)$. But $X_1$, $X_2$ are regularly homotopic if and only if $\tilde{q}_1 \equiv \tilde{q}_2 (\mathrm{mod}\, 2)$, and the corresponding immersed surfaces are regularly homotopic if and only if $\mathrm{Arf}\, \tilde{q}_1 = \mathrm{Arf}\, \tilde{q}_2$ (see [24]). □

## 7 Theta characteristics and spin structures

Theorem 7 ties the notion of spin structure with other concepts from algebraic geometry. Recall that a *theta characteristic* on a Riemann surface is a divisor $D$ such that $2D$ is the canonical divisor.

THEOREM 7. *Given a Riemann surface $M$, there are natural bijections between the following sets of objects:*
  (i) *the spin structures on $M$;*
  (ii) *the complex line bundles $S$ on $M$ satisfying $S \otimes S \cong T(M)$;*
  (iii) *the theta characteristics on $M$;*

(iv) *the classes of non-identically-zero meromorphic differentials on $M$ whose zeros and poles have even orders, under the equivalence*

$$\varphi_1 \sim \varphi_2 \quad \Longleftrightarrow \quad \varphi_1/\varphi_2 = h^2 \text{ for some meromorphic function } h \text{ on } M.$$

*Proof.* (i) $\Longleftrightarrow$ (ii): Given a line bundle $S$ on $M$ satisfying $S \otimes S \cong T(M)$, $S$ is a spin structure with mapping $\mu : S \longrightarrow S \otimes S$ defined by $\mu(s) = s \otimes s$. Conversely, given a spin structure $S$ on $M$, the map $\mu(s) \mapsto s \otimes s$ is well-defined and a vector-bundle isomorphism, so $T(M)$ is isomorphic to $S \otimes S$.

(ii) $\Longleftrightarrow$ (iii): Via the natural correspondence between the line bundles on $M$ with the divisor classes, this set of line bundles is bijective with with the theta characteristics.

(iii) $\Longleftrightarrow$ (iv): Again, there is a natural bijection between the meromorphic differentials with zeros and poles of even orders and the theta characteristics. Given such a differential $\varphi$, the corresponding theta characteristic is $\frac{1}{2}(\varphi)$. Moreover, two such differentials correspond to theta characteristics in the same linear equivalence class if and only if their ratio is the square of a meromorphic function on $M$. For

$$\varphi_1/\varphi_2 = h^2 \quad \Longleftrightarrow \quad \tfrac{1}{2}(\varphi_1) - \tfrac{1}{2}(\varphi_2) = (h).$$

$\square$

The spin structures on a compact Riemann surface are also bijective with the various translates $\vartheta\begin{bmatrix}a_0\\b_0\end{bmatrix}$ of the theta functions on the surface (see [21] for the definition of $\vartheta\begin{bmatrix}a_0\\b_0\end{bmatrix}$).

## 8 Spin structures on hyperelliptic Riemann surfaces

In the special case of a hyperelliptic Riemann surface, the spin structures and their corresponding quadratic forms are computed explicitly.

THEOREM 8. *Let*

$$M = \left\{ [x_1, x_2, x_3] \in \mathbb{CP}^2 \ \Big|\ x_2^2 x_3^{2g-1} = \prod_{i=1}^{2g+1}(x_1 - a_i x_3) \right\}$$

*be a hyperelliptic Riemann surface of genus $g$, where $A = \{a_1, \ldots, a_{2g+1}\} \subset \mathbb{C}$ is a set of $2g+1$ distinct points. Let $z = x_1/x_3$ and $w = x_2/x_3$. For each subset $B \subseteq A$, define*

$$f_B(z) = \prod_{a \in B}(z - a) \quad \text{and} \quad \eta_B = f_B(z)dz/w.$$

*Then*

(i) Any differential $\eta_B$ represents a spin structure in the sense of Theorem 7.

(ii) The set of $2^{2g}$ meromorphic differentials

$$\{\eta_B \mid B \subseteq A, \#B \leq g\}$$

represent the $2^{2g}$ distinct spin structures on $M$.

(iii) With $q$ the quadratic form corresponding to $\eta_B$, let $\gamma$ be a curve in $M$ whose projection to the $z$-plane is a Jordan curve which avoids $\infty$ and $A$, and let $C \subseteq A$ be the set of branch points which lie in the region enclosed by $\gamma$ (so $\#C$ is even). Then

$$q([\gamma]) = \#(B \cap C) + \tfrac{1}{2}\#C \pmod 2.$$

(iv) With $\eta_B$ and $q$ as in (iii),

$$\operatorname{Arf} q = \begin{cases} +1 & \text{if } 2g - 2\#B + 1 \equiv \pm 1 \pmod 8, \\ -1 & \text{if } 2g - 2\#B + 1 \equiv \pm 3 \pmod 8. \end{cases}$$

*Proof of (i).* Let $P_i = P_{a_i} = [a_i, 0, 1]$ and $P_\infty = [0, 1, 0]$ be the branch points of the two-sheeted cover $z : M \longrightarrow \mathbb{CP}^1$. Then the divisor of $\eta_B$ is

$$2\left((g - \#B - 1)P_\infty + \sum_{a \in B} P_a\right).$$

Since this divisor is even, the differential represents a spin structure by Theorem 7.

*Proof of (ii).* Note that there are $\binom{2g+1}{r}$ differentials in the $(r+1)^{\text{th}}$ row, totaling $\sum_{r=0}^{g} \binom{2g+1}{r} = 2^{2g}$. All but those in the last row are holomorphic.

In order to prove that these differentials represent distinct spin structures, we first compute the relations on the divisors of the form $\sum k_i P_i + k_\infty P_\infty$. Two such divisors are equivalent if and only if there is a meromorphic function $M$ whose divisor is their difference. Since the functions $w$ and $z - a_i$ have respective divisors

$$\begin{aligned}(w) &= P_1 + \cdots + P_{2g+1} - (2g+1)P_\infty, \\ (z - a_i) &= 2P_i - 2P_\infty,\end{aligned}$$

we have the independent relations

$$\begin{aligned}P_1 + \cdots + P_{2g+1} &\equiv (2g+1)P_\infty, \\ 2P_i &\equiv 2P_\infty \quad (i = 1, \ldots, 2g+1).\end{aligned}$$

To show that there are no other relations independent of these, let $\sum k_i P_i + k_\infty P_\infty \equiv 0$ be a relation. Then $\sum k_i = k_\infty$, and by the relations above, we may assume each $k_i$ is 0 or 1. Hence the general relation may be assumed to be of the form $D - dP_\infty \equiv 0$, where $D$ is a sum of distinct $P_i \in A$, and $d = \#D$. Let $h$ be a function with divisor $D - dP_\infty$. Since the only pole of $h$ is at $P_\infty$, $h$ is a polynomial in $z$ and $w$, so there are polynomial functions $f_1$ and $f_2$ of $z$ such that

$$h(z,w) = f_1(z) + wf_2(z).$$

Then

$$2g+1 \geq d = -\operatorname{ord}_{P_\infty} h = -\operatorname{ord}_{P_\infty}(f_1 + wf_2) \geq -\operatorname{ord}_{P_\infty} wf_2 = \deg f_2 + 2g + 1.$$

Thus $d = 2g+1$, and $D = P_1 + \cdots + P_{2g+1}$, so no new relation can exist.

We want to show that $\eta_{B_1}$ and $\eta_{B_2}$ represent identical spin structures if and only if $B_1 = B_2$ or $B_1 = B_2'$, where the prime notation $C'$ designates the complement $A \setminus C$ in $A$. If $B_1 = B_2$, then this is clear; if $B_1 = B_2'$, then $\eta_{B_2}/\eta_{B_1} = (f_2/w)^2$ is a square of a meromorphic function on $M$, and so $\eta_{B_1}$ and $\eta_{B_2}$ represent the same spin structure by Theorem 7.

Conversely, suppose that $\eta_{B_1}$ and $\eta_{B_2}$ represent the same spin structure. Then by Theorem 7, $\eta_{B_2}/\eta_{B_1} = h^2$ for some meromorpic function $h$ on $M$. But

$$2(h) = (h^2) = (\eta_{B_2}/\eta_{B_1}) = 2((d_2 - d_1)P_\infty + D_2 - D_1),$$

where $D_1 = \sum_{a \in B_1} P_a$, $D_2 = \sum_{a \in B_2} P_a$, $d_1 = \#B_1$, and $d_2 = \#B_2$. So $(d_2 - d_1)P_\infty + D_2 - D_1 \equiv 0$. By the relations (8), this divisor is equivalent to $\sum_{a \in B_1 \circ B_2} P_a - \#(B_1 \circ B_2)P_\infty$, where $B_1 \circ B_2$ is the symmetric difference $(B_1 \cup B_2) \setminus (B_1 \cap B_2)$. Since the relations (8) generate all such relations, it follows that $B_1 \circ B_2$ is either $\emptyset$ or $A$, that is that $B_1 = B_2$ or $B_1 = B_2'$.

*Proof of (iii).* It follows from the definition of $q$ that $q([\gamma])$ is the degree (mod 2) of the map $f(z)/w$ thought of as a map from the curve $\gamma$ on $M$ to $\mathbb{C} \setminus \{0\}$. Let $h = (f/w)^2$. Then

$$\deg h = \sum_{h(p)=0} \operatorname{ord}_p h + \sum_{h(p)=\infty} \operatorname{ord}_p h,$$

the sums being restricted to points within $\gamma$. This computes to

$$\deg h = \#(B \cap C) - \#(B \cup C) = 2(\#(B \cap C) - \tfrac{1}{2}\#C),$$

which shows that

$$q([\gamma]) = \#(B \cap C) + \tfrac{1}{2}\#C \pmod{2}.$$

*Proof of (iv).* In order to compute $\operatorname{Arf} q$, we first compute $\sum q(\alpha)$, where $\alpha$ ranges over $H_1(M, \mathbb{Z}_2)$. Correspondingly, the set of branch points $C$ in the region enclosed by $\alpha$ range over the subsets of $A$ of even cardinality. Hence $\sum q(\alpha)$ is the number of such subsets for which $q(\alpha) = 1$, that is, for which

$$\#(B \cap C) - \#(B' \cap C) \equiv 2 \pmod{4}.$$

The set of such subsets is

$$\{R \cup S \mid R \subseteq B, S \subseteq B', \#R - \#S \equiv 2 \pmod 4\}.$$

The cardinality of this set is

$$\sum q(\alpha) = \sum_{i-j \equiv 2} \binom{b}{i}\binom{b'}{j},$$

where $b = \#B$, $b' = \#B'$, and the sum is over $i$ and $j$ with $i - j \equiv 2 \pmod 4$.

To compute this sum, define

$$\xi(c,k) = \sum_{i \equiv k} \binom{c}{i}.$$

Then

$$\sum q(\alpha) = \sum_i \binom{b}{i} \sum_{j \equiv i+2} \binom{b'}{j} = \sum_k \binom{b}{i} \xi(b', j+2)$$

$$= \sum_{p=0}^{3} \sum_n \binom{b}{4n+p} \xi(b', p+2) = \sum_{p=0}^{3} \xi(b,p)\xi(b', p+2).$$

Using a fact about Pascal's triangle

$$\xi(c,k) = 2^{(c-2)/2}\left(2^{(c-2)/2} + \cos\tfrac{\pi}{4}(c-2k)\right),$$

we have

$$\sum q(\alpha) = 2^{(2g-3)/2} \sum_{p=0}^{3} \left(2^{(b-2)/2} + \cos\tfrac{\pi}{4}(b-2p)\right)\left(2^{(b'-2)/2} + \cos\tfrac{\pi}{4}(b'-2p)\right)$$

$$= 2^{g-1}\left(2^g - \frac{1}{\sqrt{2}}\sum_{p=0}^{3} \cos\tfrac{\pi}{4}(b-2p)\cos\tfrac{\pi}{4}(b'-2p)\right)$$

$$= 2^{g-1}\left(2^g - \sqrt{2}\cos\tfrac{\pi}{4}(2g - 2b + 1)\right)$$

$$= \begin{cases} 2^{g-1}(2^g - 1) & \text{if } 2g - 2b + 1 \equiv \pm 1 \pmod 8, \\ 2^{g-1}(2^g + 1) & \text{if } 2g - 2b + 1 \equiv \pm 3 \pmod 8. \end{cases}$$

Since $(-1)^t = 1 - 2t$ for $t = 0$ or $1$,

$$\operatorname{Arf} q = \frac{1}{2^g} \sum (-1)^{q(\alpha)} = \frac{1}{2^g}(2^{2g} - 2\sum q(\alpha))$$

is $+1$ or $-1$ according as $2g - 2b + 1$ is $\pm 1$ or $\pm 3 \pmod{8}$. □

As an example, we compute the values of $q$ for the four spin structures on a Riemann torus $T$. Let $\mathbb{C}/\{2\omega_1, 2\omega_3\} = \text{Jac}(T)$ be the Jacobian for $T$, and let $e_i = \wp(\omega_i)$, $\omega_2 = \omega_1 + \omega_3$. Then $\varphi(u) = (\wp(u), \wp'(u))$ maps the Jacobian to the Riemann surface $M$ defined by $w^2 = 4(z - e_1)(z - e_2)(z - e_3)$. The differentials as in (ii) of the above theorem pull back to

$$\begin{aligned} du &= \varphi^*(dz/w), \\ (\wp(u) - e_i)du &= \varphi^*((z - e_i)dz/w). \end{aligned}$$

With $\alpha_i$ the generator of $H_1(T, \mathbb{Z}_2)$ defined by

$$\alpha_i : [0, 1] \longrightarrow \text{Jac}(T), \quad \alpha_i(t) = 2t\omega_i,$$

the values of $q$ and $\text{Arf} \, q$ are tabulated.

Table 1: Values of $q$ and $\text{Arf} \, q$ for spin structures on tori

| $\eta$ | $q_\eta(0)$ | $q_\eta(\alpha_1)$ | $q_\eta(\alpha_2)$ | $q_\eta(\alpha_3)$ | $\text{Arf} \, q_\eta$ |
|---|---|---|---|---|---|
| $du$ | 0 | 1 | 1 | 1 | $-1$ |
| $(\wp(u) - e_1)du$ | 0 | 1 | 0 | 0 | $+1$ |
| $(\wp(u) - e_2)du$ | 0 | 0 | 1 | 0 | $+1$ |
| $(\wp(u) - e_3)du$ | 0 | 0 | 0 | 1 | $+1$ |

An immersion corresponding to $q$ for which $\text{Arf} \, q = +1$ is regularly homotopic to the torus standardly embedded in $\mathbb{R}^3$. The value $\text{Arf} \, q = -1$ corresponds to the twisted torus, which can be realized as the "diagonal" double covering of the standardly embedded torus as shown, but is not regularly homotopic to an embedding.

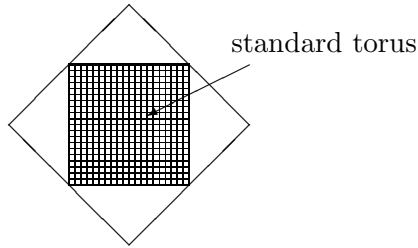

Figure 4: The twisted torus

# 9  Group action on spinors

The largest linear group acting on $Q$ is the "linear conformal group"

$$\mathbb{C}^* \times \mathrm{SO}(3,\mathbb{C}).$$

The orbit of an immersion into $Q$ under this action is an 8-real-dimensional family of immersions (which action, however, will not respect the vanishing of periods — see section 10). The subgroup

$$\mathbb{R}^+ \times \mathrm{SO}(3,\mathbb{R})$$

is the group of similarity transformations of Euclidean 3-space. Hence the homogeneous space

(9.6) $\quad (\mathbb{C}^* \times \mathrm{SO}(3,\mathbb{C})) / (\mathbb{R}^+ \times \mathrm{SO}(3,\mathbb{R})) \cong S^1 \times (\mathrm{SO}(3,\mathbb{C})/\mathrm{SO}(3,\mathbb{R}))\,.$

is the 4-real-dimensional parameter space of non-similar surfaces in the above orbit.

The $S^1$ factor gives rise to the well-known "associate family" of minimal surfaces, which are locally isometric and share a common Gauss map. The other factor has a simple (though apparently less known) geometric interpretation as well. The Gauss map is the ratio of two spinors, so $\mathrm{SO}(3,\mathbb{C}) \cong \mathrm{PSL}(2,\mathbb{C})$ acts on the Gauss map via post-composition with a fractional linear transformation of $S^2$; indeed, the quotient by $\mathrm{SO}(3,\mathbb{R}) \cong \mathrm{PSU}(2)$ leaves the hyperbolic three-space $H^3$, so the second factor can be thought of as the non-rigid Möbius deformations of the Gauss map.

The above observations are justified by the following well-known fact (see, for example, [8], [25]).

<u>Theorem 9.</u>  *There is a unique two-fold covering homomorphism*

$$T : \mathrm{GL}(2,\mathbb{C}) \longrightarrow \mathbb{C}^* \times \mathrm{SO}(3,\mathbb{C})$$

*such that for any $A \in \mathrm{GL}(2,\mathbb{C})$,*

(9.7) $\qquad\qquad\qquad\qquad T(A)\sigma = \sigma A,$

*where $\sigma : \mathbb{C}^2 \longrightarrow Q$ is as in equation (4.5), and $A$ and $T(A)$ act by left multiplication on $C^2$ and $Q$ respectively. Moreover, $T$ is a two-fold covering homomorphism when restricted to the following groups:*

$$\begin{aligned}
T &: \mathrm{GL}(2,\mathbb{C}) \longrightarrow \mathbb{C}^* \times \mathrm{SO}(3,\mathbb{C}), \\
T &: \mathrm{SL}(2,\mathbb{C}) \longrightarrow \mathrm{SO}(3,\mathbb{C}), \\
T &: \mathbb{R}^* \times \mathrm{SU}(2) \longrightarrow \mathbb{R}^+ \times \mathrm{SO}(3,\mathbb{R}), \\
T &: \mathrm{SU}(2) \longrightarrow \mathrm{SO}(3,\mathbb{R}).
\end{aligned}$$

*Proof.* We define $T$ and omit many of the calculations. $\mathbb{C}^3$ may be identified with the set $\Gamma$ of traceless $2 \times 2$ complex matrices via

$$(x_1, x_2, x_3) \longleftrightarrow \begin{pmatrix} x_3 & -x_1 + ix_2 \\ -x_1 - ix_2 & -x_3 \end{pmatrix} = X,$$

with the subset $\mathbb{R}^3 \subset \mathbb{C}^3$ identified with $\Gamma_\mathbb{R} = \{X \in \Gamma \mid X = \overline{X}^t\}$ The inner product on $\mathbb{C}^3$ becomes

$$\langle X, Y \rangle = \sum_1^3 x_i y_i = \tfrac{1}{2} \operatorname{tr} XY,$$

and

$$\langle X, X \rangle = \frac{1}{2} \operatorname{tr} X^2 = -\det X,$$

so $Q \subset \mathbb{C}^3$ is identified with

$$\Gamma_Q = \{X \in \Gamma \mid \det X = 0\}.$$

Similarly, $\mathbb{C}^2$ may be identified with the set $\Delta$ of matrices of the form

$$\begin{pmatrix} x_1 & x_1 \\ x_2 & x_2 \end{pmatrix}.$$

The map $\sigma : \mathbb{C}^2 \longrightarrow Q$ becomes under these identifications $\sigma : \Delta \longrightarrow \Gamma_Q$ given by $\sigma(X) = XKX'$, where $K = \begin{pmatrix} 0 & -1 \\ 1 & 0 \end{pmatrix}$, and $X'$ denotes the classical adjoint

$$\begin{pmatrix} a & b \\ c & d \end{pmatrix}' = \begin{pmatrix} d & -b \\ -c & a \end{pmatrix}$$

satisfying $XX' = X'X = (\det X)I$ and $(XY)' = Y'X'$.

Then in order to satisfy equation (9.7), $T$ must be defined, for $X \in \Gamma$, by

$$T(A)X = AXA'.$$

It follows that $T(A)$ is linear and maps $\Gamma$ to itself, and that $T : \operatorname{GL}(2, \mathbb{C}) \longrightarrow \operatorname{GL}(3, \mathbb{C})$ is a homomorphism with kernel $\{\pm I\}$. That $T$ maps into the four groups listed follows from the equation

$$\langle T(A)X, T(A)Y \rangle = (\det A)^2 \langle X, Y \rangle$$

and the fact that $T(A)(\Gamma_\mathbb{R}) = \Gamma_\mathbb{R}$ for $A \in \mathbb{R}^* \times \operatorname{SU}(2)$. □

Lifting the group action on $Q$ to an action on $\mathbb{C}^2 \setminus \{0\}$ via $T$, the homogeneous space in equation (9.6) can also be written

(9.8)  $(\operatorname{GL}(2, \mathbb{C})) / (\mathbb{R}^* \times \operatorname{SU}(2)) \cong S^1 \times (\operatorname{SL}(2, \mathbb{C})/\operatorname{SU}(2)) \cong S^1 \times H^3,$

where $H^3$ is hyperbolic three-space.

## 10 Periods

Given an immersion $X : M \longrightarrow \mathbb{R}^3$, the *period* around a simple closed curve $\gamma \subset M$ is the vector in $\mathbb{C}^3$
$$\int_\gamma \partial X.$$
If the real part of a period is not $(0,0,0)$, the resulting surface is periodic and does not have finite total curvature. It is a considerable problem to "kill the periods" — that is, choose parameters so that the integrals around every simple closed curve in $M$ generates purely imaginary period vectors. Non-zero periods can arise along two kinds of simple closed curves:
- a simple closed curve around an end $p \in M$,
- a non-trivial simple closed curve in $H_1(M, \mathbb{Z})$.

For a simple closed curve $\gamma$ around an end $p \in M$,
$$\int_\gamma \partial X = 2\pi i \operatorname{res}_p \partial X.$$
This integral is zero at embedded planar ends.

Using the spinor representation, the condition that a period along a closed curve $\gamma \subset M$ be pure imaginary can be expressed by
$$\int_\gamma \left(s_1^2 - s_2^2, i(s_1^2 + s_2^2), 2s_1 s_2\right) \in i\mathbb{R}^3,$$
equivalent to

(10.9)
$$\begin{aligned} \int_\gamma s_1^2 &= \overline{\int_\gamma s_2^2} \\ \int_\gamma s_1 s_2 &\in i\mathbb{R}. \end{aligned}$$

These equations are preserved by the group $\mathbb{R}^* \times \operatorname{SU}(2)$ of homotheties.

## 11 Spin structures and nonorientable surfaces

To deal with immersions of a nonorientable manifold $M$ into space, we pass to the oriented two-fold cover of $M$. The following rather technical results are required in Part III. Without proof we state:

LEMMA 10. *Let*

$A : S^2 \longrightarrow S^2$ be the antipodal map,
$A_* : T(S^2) \longrightarrow T(S^2)$ the derivative of $A$,
$\hat{A}_* : \operatorname{Spin}(S^2) \longrightarrow \operatorname{Spin}(S^2)$ one of the lifts of $A_*$ to $\operatorname{Spin}(S^2)$ .

Then in the coordinates of section 5 we have

$$\begin{aligned} A_* &= \text{Conj,} \\ \hat{A}_* &= \pm \begin{pmatrix} 0 & i \\ -i & 0 \end{pmatrix} \circ \text{Conj.} \end{aligned}$$

The lifts of the antipodal maps are shown in Figure 5.

$$\begin{array}{ccc} \text{Spin}(S^2) & \xrightarrow{\hat{A}_*} & \text{Spin}(S^2) \\ \downarrow & & \downarrow \sigma \\ T(S^2) & \xrightarrow{A_*} & T(S^2) \\ \downarrow & & \downarrow \\ S^2 & \xrightarrow{A} & S^2 \end{array}$$

Figure 5: Lifts of the antipodal map

<u>THEOREM 11.</u> *Let $M$ be a nonorientable Riemann surface, and $X : M \longrightarrow \mathbb{R}^3$ a conformal minimal immersion of $M$ into space.*

*Let $\pi : \widetilde{M} \longrightarrow M$ be an oriented double cover of $M$, and $\widetilde{X} = X \circ \pi$ the lift of $X$ to this cover. Let $I : \widetilde{M} \longrightarrow \widetilde{M}$ the order-two deck transformation for the cover. With $\omega = \partial \widetilde{X}$, and in the notation of Lemma 10, we have*
  *(i) $gI = Ag$,*
  *(ii) $\omega I_* = A_* \omega$,*
  *(iii) $\hat{\omega} \hat{I}_* = \pm \hat{A}_* \hat{\omega}$.*

*Proof.* Since $\widetilde{X}$ is the lift of $X$, we have that $\widetilde{X} = \widetilde{X} I$. Hence $\text{Re } \omega = \text{Re } \partial \widetilde{X} = d\widetilde{X} = d\widetilde{X} I_* = \text{Re } \omega I_*$. Since $X$ is a conformal minimal immersion, $\omega$ is holomorphic. Hence $\omega$ and $\omega I_*$ are either equal or conjugate. But $I_*$ is orientation reversing, so they are conjugate and $\omega I_* = \overline{\omega} = A_* \omega$, proving (ii).

From
$$gI\pi_1 = g\pi_1 I_* = \pi_3 \omega I_* = \pi_3 A_* \omega = A\pi_3 \omega = Ag\pi_1$$
and the surjectivity of $\pi_1$, (ii) follows. Similarly, from
$$\sigma \hat{\omega} \hat{I}_* = \omega \pi_2 \hat{I}_* = \omega I_* \pi_2 = A_* \omega \pi_2 = A_* \sigma \hat{\omega} = \sigma \hat{A}_* \hat{\omega}$$

(iii) follows. □

We remark that the proper treatment of nonorientable surfaces should really be via "pin" structures (Pin($n$) being the corresponding two-sheeted covering group of O($n$)), and that in this case we should have an analytic formula for the full $\mathbb{Z}_8$-valued Arf invariant of the associated $\mathbb{Z}_4$-valued quadratic form on $H^1(M, \mathbb{Z}_2)$.

# Part II
## Minimal Immersions with Embedded Planar Ends

## 12  Embedded planar ends

The first section of this part of our paper discusses the behavior of a minimal immersion at an embedded planar end. Lemma 12 translates this geometric behavior to a necessary and sufficient algebraic condition on the order and residue of the immersion at the end. Arising naturally from this algebraic condition is a certain vector subspace $\mathcal{K}$ of holomorphic spin-vector fields (sections of a spin structure) which generates all minimal surfaces with embedded planar ends (Theorem 15). More precisely, two sections chosen from $\mathcal{K}$ form the spinor representation of a minimal surface, and conversely, any such surface must arise this way. However, such a surface is usually periodic, and possibly a branched immersion. In order to compute $\mathcal{K}$ explicitly, a skew-symmetric bilinear form $\Omega$ is next defined (Definition 16) whose kernel is closely related to the space $\mathcal{K}$. In Part III, this form is used to prove existence and non-existence theorems for a variety of examples.

## 13  Algebraic characterization of embedded planar ends

The geometric condition that an end of a minimal immersion be embedded and planar can be translated to algebraic conditions (see for example [4]). Let $X : D \backslash \{p\} \longrightarrow \mathbb{R}^3$ be a conformal minimal immersion of an open disk $D$ punctured at $p$ such that $\lim_{q \to p} |X(q)| = \infty$. The image under $X$ of a small neighborhood of $p$ (and by association, $p$ itself) is what we shall refer to as an *end*. The behavior of the end is determined by the residues and the orders of the poles of $\partial X$ at $p$ as follows.

Let $\psi_1$, $\psi_2$, $\psi_3$ be defined by

$$\partial X = (\psi_1 - \psi_2, i(\psi_1 + \psi_2), 2\psi_3).$$

The Gauss map for this immersion (see [22]) is

$$g = \psi_3/\psi_1 = \psi_2/\psi_3.$$

First note that for $X$ to be well-defined, we must have for any closed curve $\gamma$, which winds $k$ times around $p$,

$$0 = \mathrm{Re} \int_\gamma \partial X = k \, \mathrm{Re}\, (2\pi i \operatorname{res}_p \partial X),$$

and so $\text{res}_p \partial X$ must be *real*. Assume this, and assume initially that the limiting normal to the end is upward (that is $g(p) = \infty$). In this case,

$$\text{ord}_p \psi_2 < \text{ord}_p \psi_3 < \text{ord}_p \psi_1,$$

so the first two coordinates of $X(q)$ grow faster than does the third as $q \to p$.

It follows that $\text{ord}_p \psi_2$ cannot be $-1$, because if it were then

$$\text{res}_p \partial X = (-\text{res}_p \psi_2, i\,\text{res}_p \psi_2, 0)$$

would not be real. Hence $\text{ord}_p \psi_2 \leq -2$. The image under $X$ of a small closed curve around $p$ is a large curve which winds around the end $|\text{ord}_p \psi_2| - 1$ times. The end is embedded precisely when $\text{ord}_p \psi_2 = -2$.

If an end is embedded, its behavior is determined by the vanishing or non-vanishing of the residues of $\partial X$. For an embedded end, $-2 = \text{ord}_p \psi_2 < \text{ord}_p \psi_3$, so $\psi_3$ has either a simple pole or no pole. If $\psi_3$ has a simple pole (and hence also a residue), the end grows logarithmically relative to its horizontal radius and is a *catenoid* end. If $\psi_3$ has no pole, the end is asymptotic to a horizontal plane and is called a *planar* end. Moreover, in this latter case, $\text{res}_p \psi_2$ must vanish (again, if it did not, $\text{res}_p \partial X$ would not be real), and so $\text{res}_p \partial X = (0, 0, 0)$.

For an end in general position the same conclusions hold, because a real rotation affects neither $\text{ord}_p \partial X$ nor the reality or vanishing of $\text{res}_p \partial X$. In summary, we have

LEMMA 12. *Let $X: D\backslash\{p\} \longrightarrow \mathbb{R}^3$ be a conformal minimal immersion of a punctured disk. Then $p$ is an embedded planar end if and only if*

$$\text{ord}_p \partial X = -2 \quad \text{and} \quad \text{res}_p \partial X = 0,$$

*where $\text{ord}_p \partial X$ denotes the minimum order at $p$ of the three coordinates of $\partial X$.*

## 14 Embedded planar ends and spinors

The conditions in the lemma above can be translated into conditions on the spinor representation of the minimal immersion. This leads to the definition of a space $\mathcal{K}$ of spin-vector fields, pairs of which form the spinor representation satisfying the required conditions.

Throughout the rest of Part II, the following notation is fixed:

(14.10)
$M$ is a compact Riemann surface,
$S$ is a spin structure on $M$,
$P = [p_1] + \ldots + [p_n]$ is a divisor of $n$ distinct points.

The points $p_1, \ldots, p_n$ will eventually be the ends of a minimal immersion of $M$ whose spinor representation will be a pair of sections of $S$.

Let $H^0(M, \mathcal{O}(S))$ and $H^0(M, \mathcal{M}(S))$ denote respectively the vector spaces of holomorphic and meromorphic sections of $S$. Define

(14.11)
$$\begin{aligned}
\mathcal{F} &= \mathcal{F}_{M,S,P} = \{s \in H^0(M, \mathcal{M}(S)) \mid (s) \geq -P\} \\
\mathcal{H} &= \mathcal{H}_{M,S} = H^0(M, \mathcal{O}(S)) \\
\mathcal{K} &= \mathcal{K}_{M,S,P} = \{s \in \mathcal{F} \mid \operatorname{ord}_p s \neq 0 \text{ and } \operatorname{res}_p s^2 = 0 \text{ for all } p \in \operatorname{supp} P\}.
\end{aligned}$$

THEOREM 13.   *Let $X : M \longrightarrow \mathbb{R}^3$ be a minimal immersion with spinor representation $(s_1, s_2)$. Then $p \in M$ is an embedded planar end if and only if $s_1, s_2 \in \mathcal{K}$ and at least one of $s_1$, $s_2$ has a pole at $p$.*

*Proof.* By Lemma 12, $p$ is an embedded planar end if and only if

$$\operatorname{ord}_p \partial X = -2 \quad \text{and} \quad \operatorname{res}_p \partial X = 0.$$

The first of these equations is equivalent to the condition

$$s_1, s_2 \in \mathcal{F}, \text{ and at least one of } s_1, s_2 \text{ has a pole at } p.$$

Given this, the conditions $\operatorname{ord}_p s_1 \neq 0$, $\operatorname{ord}_p s_2 \neq 0$ in the definition of $\mathcal{K}$ follow because if one were 0, the other would be $-1$, giving $s_1 s_2$ a nonvanishing residue.

It remains only to show that for $s_1, s_2 \in \mathcal{F}$,

$$\operatorname{res}_p s_1^2 = 0, \ \operatorname{res}_p s_2^2 = 0 \quad \Longrightarrow \quad \operatorname{res}_p s_1 s_2 = 0.$$

This is an application of the following lemma.                                          □

LEMMA 14.   *Let $S$ be a spin structure on a closed Riemann surface $M$, and let $s_1$, $s_2$ be meromorphic sections of $S$ with $\operatorname{ord}_p s_1 \geq -1$, $\operatorname{ord}_p s_2 \geq -1$ for some $p \in M$. Then*

$$2 \operatorname{res}_p s_1 s_2 = \begin{cases} \left[\frac{s_2}{s_1}\right]_p \operatorname{res}_p s_1^2 + \left[\frac{s_1}{s_2}\right]_p \operatorname{res}_p s_2^2 & \text{if } \operatorname{ord}_p s_1 = \operatorname{ord}_p s_2, \\ 0 & \text{if } |\operatorname{ord}_p s_1 - \operatorname{ord}_p s_2| \geq 2. \end{cases}$$

*Proof.* With $z$ a complex coordinate near $p$ satisfying $z(p) = 0$, let $\varphi$ be a section of $S$ satisfying $\varphi^2 = dz$. Also let

$$s_1 = \left(\frac{a_{-1}}{z} + a_0 + \ldots\right) \varphi,$$

$$s_2 = \left(\frac{b_{-1}}{z} + b_0 + \ldots\right) \varphi$$

be the expansions of $s_1$ and $s_2$. Then

$$\begin{aligned}
\operatorname{res}_p s_1^2 &= 2a_{-1}a_0, \\
\operatorname{res}_p s_2^2 &= 2b_{-1}b_0, \\
\operatorname{res}_p s_1 s_2 &= a_{-1}b_0 + a_0 b_{-1}.
\end{aligned}$$

In case $\operatorname{ord}_p s_1 = \operatorname{ord}_p s_2$, then either $a_{-1} \neq 0$, $b_{-1} \neq 0$, so that

$$\begin{aligned}
\left[\frac{s_2}{s_1}\right]_p \operatorname{res}_p s_1^2 + \left[\frac{s_1}{s_2}\right]_p \operatorname{res}_p s_2^2 &= \frac{b_{-1}}{a_{-1}}(2a_{-1}a_0) + \frac{a_{-1}}{b_{-1}}(2b_{-1}b_0) \\
&= 2(b_{-1}a_0 + a_{-1}b_0) \\
&= 2\operatorname{res}_p s_1 s_2,
\end{aligned}$$

or $a_{-1} = b_{-1} = 0$, and the three residues vanish.

In case $|\operatorname{ord}_p s_1 - \operatorname{ord}_p s_2| \geq 2$, then $a_{-1} = a_0 = 0$ or $b_{-1} = b_0 = 0$, so again the three residues vanish. So in each case the formula is verified. □

## 15  $\mathcal{F}$ and $\mathcal{K}$ as vector spaces

The following theorem develops some of the properties of the spaces $\mathcal{F}$ and $\mathcal{K}$. The most important of these is that $\mathcal{K}$ is in fact a vector space. In this section we write $K$ for the holomorphic cotangent bundle (that is, the *canonical* line bundle) of $M$.

<u>Theorem 15.</u>  *With $M$, $P$, and $S$ as in equation (14.10), and $\mathcal{F}$, $\mathcal{H}$, and $\mathcal{K}$ as in equation (14.11), let $g = \operatorname{genus}(M)$ and $n = \deg P$. Then*
  (i) *if $n \geq g$, then $\dim \mathcal{F} = n$;*
 (ii) *$\mathcal{K}$ and $\mathcal{H}$ are subspaces of $\mathcal{F}$;*
(iii) *if $n \geq g$, then $\mathcal{K} \cap \mathcal{H} = 0$.*

*Proof of (i):* The dimension of $\mathcal{F}$ can be computed by means of the Riemann-Roch theorem (see, for example, [10]) which states

$$\dim H^0(M, L) - \dim H^0(M, K \otimes L^*) = \deg L - g + 1$$

for an arbitrary line bundle $L$. Let $R$ be the line bundle corresponding to the divisor $P$, and let $L = S \otimes R$. Then:
- $H^0(M, L) \cong \mathcal{F}$ by the isomorphism $s \otimes r \mapsto s$, where $r$ is a section of $R$ with divisor $P$;
- $H^0(M, K \otimes L^*) = 0$, since $\deg(K \otimes L^*) = g - 1 - n$, which is negative by hypothesis;

- $\deg L - g + 1 = n$;

from which it follows that

$$\dim \mathcal{F} = \dim H^0(M, L) = n.$$

*Proof of (ii):* $\mathcal{H} \subseteq \mathcal{F}$ is a subspace by definition. To show that $\mathcal{K} \subseteq \mathcal{F}$ is a subspace, let $s \in \mathcal{K}$ and $p \in P$, let $z$ be a conformal coordinate near $p$ with $z(p) = 0$, and let $\varphi$ be a section of $S$ which satisfies $\varphi^2 = dz$. Let

$$s = \left(\frac{a_{-1}}{z} + a_0 + \ldots\right)\varphi$$

be the expansion of $s$. Then

$$s^2 = \left(\frac{a_{-1}^2}{z^2} + \frac{2a_{-1}a_0}{z} + \ldots\right) dz,$$

so that

$$\operatorname{res}_p s^2 = 0 \quad \Longleftrightarrow \quad a_{-1} = 0 \text{ or } a_0 = 0$$

and

$$\operatorname{ord}_p s^2 \neq 0 \quad \Longleftrightarrow \quad a_{-1} \neq 0 \text{ or } a_0 = 0.$$

Together, these two conditions are equivalent to the condition

$$a_0 = 0.$$

Thus

(15.11) $\qquad s \in \mathcal{K} \quad \Longleftrightarrow \quad$ the constant term in the expansion of $s$ vanishes at each $p \in P$.

$\mathcal{K}$ is a vector space because if $s_1$, $s_2$ satisfy condition (15.11), then so does any $\mathbb{C}$-linear combination of $s_1$ and $s_2$.

*Proof of (iii):* Let $s \in \mathcal{K} \cap \mathcal{H}$ be a section which is not identically zero. Since $s \in \mathcal{K}$, we have that $\operatorname{ord}_p s \neq 0$ for all $p \in \operatorname{supp} P$ — that is, at each such $p$, $s$ has either a pole or a zero. But since $s \in \mathcal{H}$, $s$ cannot have a pole at $p$, and hence has a zero, so $(s) \geq P$. Conversely, if $(s) \geq P$, then $s \in \mathcal{K} \cap \mathcal{H}$, so

$$\mathcal{K} \cap \mathcal{H} = \{s \in \mathcal{F} \mid (s) \geq P\}.$$

Thus for $s \in \mathcal{K} \cap \mathcal{H}$ not identically zero,

$$n \leq \deg s = g - 1.$$

Hence if $n \geq g$, then $\mathcal{K} \cap \mathcal{H} = 0$. □

## 16 A bilinear form $\Omega$ which kills $\mathcal{K}$

In order to understand the vector space $\mathcal{K}$ more explicitly, a skew-symmetric bilinear form $\Omega$ is defined whose kernel contains $\mathcal{K}$. This form may then be used in many cases to compute $\mathcal{K}$, and thereby moduli spaces of minimal surfaces with embedded planar ends.

DEFINITION 16. *With $M$, $P$, and $S$ as in equation (14.10) define $\Omega = \Omega_{M,P,S} : \mathcal{F} \times \mathcal{F} \longrightarrow \mathbb{C}$ by*

$$\Omega(s_1, s_2) = \sum_{p \in P} \xi(p; s_1, s_2),$$

where

$$\xi(p; s_1, s_2) = \begin{cases} \dfrac{1}{2}\left[\dfrac{s_2}{s_1}\right]_p \mathrm{res}_p s_1^2 & \text{if } \dfrac{s_2}{s_1} \text{ does not have a pole at } p, \\ \mathrm{res}_p s_1 s_2 & \text{if } \dfrac{s_2}{s_1} \text{ has a pole at } p. \end{cases}$$

The form $\Omega$ can be computed as follows: for $p \in P$, let $z$ be a conformal coordinate near $p$ with $z(p) = 0$, let $\varphi^2 = dz$, and let

$$s_1 = \left(\frac{a_{-1}}{z} + a_0 + \dots\right)\varphi,$$
$$s_2 = \left(\frac{b_{-1}}{z} + b_0 + \dots\right)\varphi$$

be the expansions of $s_1$ and $s_2$. Then

(16.12) $$\xi(p; s_1, s_2) = b_{-1} a_0.$$

THEOREM 17. *With $\mathcal{H}$, $\mathcal{K}$ as in equation (14.11), $\Omega$ satisfies the following:*
 (i) *$\Omega$ is a skew-symmetric bilinear form on $\mathcal{F}$;*
 (ii) *$\ker \Omega \supseteq \mathcal{K} + \mathcal{H}$;*
 (iii) *if $\mathcal{K} \cap \mathcal{H} = 0$, then $\ker \Omega = \mathcal{K} \oplus \mathcal{H}$;*
 (iv) *if $n = \deg P \geq \mathrm{genus}(M)$, then $\ker \Omega = \mathcal{K} \oplus \mathcal{H}$.*

*Proof of (i):* For a given $s_1, s_2 \in \mathcal{F}$, let

$$P_0 = \{\text{zeros of } s_2/s_1\} \cap P,$$
$$P_\infty = \{\text{poles of } s_2/s_1\} \cap P,$$
$$P_1 = P \setminus (P_0 \cup P_\infty),$$

so that $P_0$, $P_\infty$, $P_1$ are disjoint and their union is $P$. Then

$$\text{(16.13)} \qquad \Omega(s_1, s_2) = \frac{1}{2} \sum_{p \in P_1} \left[\frac{s_2}{s_1}\right]_p \operatorname{res}_p s_1^2 + \sum_{p \in P_\infty} \operatorname{res}_p s_1 s_2,$$

$$\Omega(s_2, s_1) = \frac{1}{2} \sum_{p \in P_1} \left[\frac{s_1}{s_2}\right]_p \operatorname{res}_p s_2^2 + \sum_{p \in P_0} \operatorname{res}_p s_1 s_2.$$

Adding the above two equations, and using the lemma in section 14 yields

$$\Omega(s_1, s_2) + \Omega(s_2, s_1) = 2 \sum_{p \in P} \operatorname{res}_p s_1 s_2,$$

which vanishes by the residue theorem, since all poles of $s_1 s_2$ are in $P$.

*Proof of (ii):* To show that that $\mathcal{K} \subseteq \ker \Omega$, let $s_1 \in \mathcal{K}$, so that $\operatorname{res}_p s_1^2 = 0$ and $\operatorname{ord}_p s_1 \neq 0$ for all $p \in P$. Let $s_2 \in \mathcal{F}$ be arbitrary, so that $\operatorname{ord}_p s_2 \geq -1$. Referring to equation (16.13), the first sum is zero because $\operatorname{res}_p s_1^2 = 0$ at each $p \in P$ by hypothesis. To show that each term in the second sum is zero, let $p \in P_\infty$ so that $\operatorname{ord}_p s_1/s_2 \geq 1$. Then

$$\operatorname{ord}_p s_1 = \operatorname{ord}_p s_2 + \operatorname{ord}_p s_1/s_2 \geq 0.$$

But $\operatorname{ord}_p s_1 \neq 0$, so $\operatorname{ord}_p s_1 \geq 1$. Then

$$\operatorname{ord}_p s_1 s_2 = \operatorname{ord}_p s_1 + \operatorname{ord}_p s_2 \geq 0,$$

so $\operatorname{res}_p s_1 s_2 = 0$.

To show that $\mathcal{H} \subseteq \ker \Omega$, let $s_1 \in \mathcal{H}$, so that $s_1^2$ has no poles, and let $s_2 \in \mathcal{F}$ be arbitrary. Then the first sum in equation (16.13) is zero because $s_1^2$ has no poles. The second sum is zero by the residue theorem — to show that all poles of $s_1 s_2$ are in $P_\infty$, note that

$$\operatorname{ord}_p s_1 s_2 = \operatorname{ord}_p s_1 + \operatorname{ord}_p s_2 \geq -\operatorname{ord}_p s_1 + \operatorname{ord}_p s_2 = \operatorname{ord}_p s_2/s_1.$$

So if $s_1 s_2$ has a pole at $p$, then so does $s_2/s_1$; this puts $p \in P_\infty$.

*Proof of (iii):* $\Omega$ can be "factored" as the composition of two maps as follows: near each point in $\operatorname{supp} P = \{p_1, \ldots, p_n\}$, choose a conformal coordinate $z_i$ with $z_i(p_i) = 0$. Let $\varphi$ be a section of $S$ satisfying $\varphi^2 = dz$, and let

$$s = \left(\frac{b_i}{z_i} + a_i + \ldots\right) \varphi$$

be the expansion of $s$ at $p_i$. Define maps $A$, $B : \mathcal{F} \longrightarrow \mathbb{C}^n$ by

$$A(s) = (a_1, \ldots, a_n),$$
$$B(s) = (b_1, \ldots, b_n).$$

By the local computation of equation (16.12), $\Omega = B^*A = -A^*B$, or with a choice of basis for $\mathcal{F}$, $\Omega = B^tA = -A^tB$ as matrices. (Note that while $\Omega$ is independent of the choice of coordinates, $A$ and $B$ are not.) Moreover, we have

$$\ker A = \mathcal{K} \quad \text{and} \quad \ker B = \mathcal{H};$$

the first is by equation (15.11), and the second is immediate from the definition of $\mathcal{H}$. (This incidentally provides another proof of (ii).)

Now let $\widehat{A} = A|_{\ker B^*A}$ and note that

$$\ker \widehat{A} = (\ker A \cap \ker B^*A) = \ker A$$

and

$$\text{image } \widehat{A} = A(\ker B^*A) = \text{image } A \cap \ker B^*.$$

Applying the rank-nullity theorem to $\widehat{A}$,

$$\begin{aligned} \dim \ker \Omega &= \dim \ker B^*A \\ &= \dim \ker \widehat{A} + \dim \text{image } \widehat{A} \\ &= \dim \ker A + \dim(\text{image } A \cap \ker B^*) \\ &\leq \dim \ker A + \dim \ker B. \end{aligned}$$

So under the assumption that $\mathcal{K} \cap \mathcal{H} = 0$,

$$\dim \ker \Omega \leq \dim \ker A + \dim \ker B = \dim \mathcal{K} + \dim \mathcal{H} = \dim \mathcal{K} \oplus \mathcal{H}.$$

But $\ker \Omega \supseteq \mathcal{K} \oplus \mathcal{H}$, so $\ker \Omega = \mathcal{K} \oplus \mathcal{H}$, proving (iii).

Part (iv) follows directly from (iii) above and Theorem 15(iii). □

# Part III
**Classification and Examples**

## 17  Genus zero

In the first half of Part III, the skew-symmetric form $\Omega$ developed in Part II is used to investigate minimal genus zero surfaces with embedded planar ends. The first two sections demonstrate the non-existence of examples with 2, 3, 5, or 7 ends, and the dimension of the moduli space of examples with 4, 6, 8, 10, 12 and 14 ends is computed. The following three sections compute explicitly the moduli spaces for the families with 4 and 6 ends, and in section 23, the moduli space of the three-ended projective planes is investigated. The remaining sections (following the heading Genus one) of Part III are devoted to constructing minimal tori and Klein bottles. All of these surfaces are found (or shown not to exist) by the following general method: after computing $\Omega$ on a simple basis, its pfaffian, which is a function of the ends, is set to zero. The resulting condition on the placement of the ends — that is, the determinantal variety — together with further conditions arising from the demand that the immersion have no periods and no branch points, forms a set of equations whose simultaneous solution (or impossibility of solution) gives the desired result.

## 18  Existence and non-existence of genus zero surfaces

The non-existence of genus zero minimal unbranched immersions with 3, 5 or 7 embedded planar ends was first proved in a case-by-case manner in [3]. The following is a new proof, using the ideas of Section 12.

THEOREM 18.  *There are no complete minimal branched or unbranched immersions of a punctured sphere into space with finite total curvature and* 2, 3, 5, *or* 7 *embedded planar ends. There exist unbranched examples with* 4, 6, *and any* $n \geq 8$ *ends.*

*Proof.* Examples with $2p$ ends ($p \geq 2$) are given in [14], and with $2p+1$ ends ($p \geq 4$) in [23]. For the cases $n = 3$, 5, or 7, by the lemma below, $2 \leq \dim \mathcal{K} \leq [\sqrt{n}] \leq 2$ (here $[q]$ denotes the greatest integer less than or equal to $q$), so $\dim \mathcal{K} = 2$, which contradicts the other statement of the lemma that $n - \dim \mathcal{K}$ is even. The case $n = 2$ is proved in [14] (or is proved likewise by the lemma).  □

   We remark that there is also a simple topological proof of the non-existence of genus zero examples with 3 ends, using ideas in [13] and [15]. The trick is to exploit the SO(3, $\mathbb{C}$)-action discussed in section 9 to deform the Gauss map — on a punctured sphere with planar ends there is no period obstruction to doing this — so that the compactified $S^2$ is transversally immersed with a unique triple-point, which is impossible. (By carefully treating the periods introduced by this explicit SO(3, $\mathbb{C}$) deformation of the Gauss map, the same kind of argument should generalize to exclude orientable minimal surfaces of any genus with three embedded planar ends — see section 27 for a proof in the case of tori.)

LEMMA 19. *Let $P$ be a divisor on the Riemann sphere $S^2$ as in equation (14.10) with $n = \deg P \geq 2$, and let $\mathcal{K} = \mathcal{K}_{S^2,S,P}$ be as in equation (14.11). Then*
  (i) *$n - \dim \mathcal{K}$ is even;*
  (ii) *If there exists a complete branched or unbranched minimal immersion of $S^2$ into space with finite total curvature and $n$ embedded planar ends in $\mathrm{supp}\,(P)$, then $2 \leq \dim \mathcal{K} \leq \sqrt{n}$.*

*Proof of (i):* By Theorem 17, $\ker \Omega = \mathcal{K} \oplus \mathcal{H}$. But $\mathcal{H} = 0$ because there are no holomorphic differentials on the sphere, so $\ker \Omega = \mathcal{K}$. Since $\Omega$ is skew-symmetric, $\mathrm{rank}\,\Omega = n - \dim \mathcal{K}$ is even (see Appendix A).

*Proof of (ii):* The sections $s_1$ and $s_2$ in the spinor representation $(s_1, s_2)$ of such a surface are independent, showing the inequality $2 \leq \dim \mathcal{K}$. To show the other inequality, let $z$ be the standard conformal coordinate on $S^2 = \mathbb{C} \cup \{\infty\}$, and let $P = \sum [a_i]$ (where the $a_i \in \mathbb{C}$ are distinct) be the divisor of the $n$ ends. Let $g_1 \eta, \ldots, g_m \eta$ be a basis for $\mathcal{K}$, where $\eta^2 = dz$. Define $f : S^2 = \mathbb{CP}^1 \longrightarrow \mathbb{CP}^{m-1}$ by

$$f = (g_1, \ldots, g_m).$$

Then $f$ is well-defined and holomorphic even at the common zeros and the common poles of $g_1, \ldots, g_m$. Let

$$h(z) = \prod (z - a_i).$$

It follows from

$$(hg_i) = (h) + (\eta) + (g_i\eta) \geq (P - n[\infty]) + [\infty] - P = -(n-1)[\infty]$$

that

$$d_0 = \deg f \leq n - 1.$$

To show that $f$ ramifies at each $a \in \mathrm{supp}\,P$, let $h_i(z) = (z - a)g_i(z)$. Then $h_i$ does not have a pole at $a$. Moreover, since by hypothesis there exists a minimal surface with ends at $\mathrm{supp}\,P$, at least one of the $g_i$ has a pole at $a$, so the $h_i$ cannot all be zero at $a$. Hence the appropriate condition that $f$ ramify at $a$ is

$$\left(h_i h_j' - h_i' h_j\right)\Big|_a = 0 \text{ for all } i, j.$$

This is satisfied because of the condition (15.11) defining $\mathcal{K}$: the expansion of $g_i$ at $a$ is

$$g_i = \frac{c_i}{z-a} + \mathrm{o}(z-a),$$

so the expansion of $h_i$ at $a$ is

$$h_i = c_i + \mathrm{o}(z-a)^2,$$

and so $h'_i(a) = 0$ for all i. Since $f$ ramifies at each $a \in \operatorname{supp} P$, we have

$$r_0 = \text{ ramification index of } f \geq n.$$

Now let $f_k : \mathbb{CP}^1 \longrightarrow \mathbb{P}(\Lambda^{k+1}\mathbb{C}^m)$ defined by $f_k = f \wedge f' \wedge \ldots \wedge f^{(k)}$ in $\mathbb{C}^m$ be the $k^{\mathrm{th}}$ *associated curve* of $f$, and use the Plücker formulas (an extension of the Riemann-Hurwitz formula — see [9]) which on $\mathbb{CP}^1$ are

$$-d_{k-1} + 2d_k - d_{k+1} - 2 = r_k,$$

where $d_k$ is the degree of $f_k$, and $r_k$ is the ramification index of $f_k$. In the table below, multiplying the numbers on the left by the inequalities on the right and adding yields

$$d_0 \geq (m+n)(m-1)/m.$$

But $n - 1 \geq d_0$, so it follows that $n \geq m^2$. □

Table 2: Values for the Plücker formulas

| | | | | | | |
|---|---|---|---|---|---|---|
| $m-1$ | | $2d_0$ | $-$ | $d_1$ | $-2 = r_0$ | $\geq n$ |
| $m-2$ | $-d_0$ | $+\ 2d_1$ | $-$ | $d_2$ | $-2 = r_1$ | $\geq 0$ |
| $\vdots$ | | $\vdots$ | | | $\vdots$ | $\vdots$ |
| 2 | $-d_{m-4}$ | $+\ 2d_{m-3}$ | $-$ | $d_{m-2}$ | $-2 = r_{m-3}$ | $\geq 0$ |
| 1 | $-d_{m-3}$ | $+\ 2d_{m-2}$ | | | $-2 = r_{m-2}$ | $\geq 0$ |

## 19  Moduli spaces of genus zero minimal surfaces

The following two theorems deal with the moduli spaces of genus zero examples.

<u>Theorem 20.</u> *Let $P$ be a divisor on $S^2$ as in equation (14.10) and $\mathcal{K} = \mathcal{K}_{S^2,S,P}$ as in equation (14.11) with $m = \dim \mathcal{K} \geq 2$. Then the space of complete minimal* branched immersions of $S^2$ into $\mathbb{R}^3$ with finite total curvature and embedded planar ends at supp $P$ is the complex $2(m-1)$-dimensional space $\mathrm{Gr}_{m,2}(\mathbb{C}) \times (S^1 \times H^3)$.

*Proof.* Each point of the Grassmanian $\mathrm{Gr}_{m,2}(\mathbb{C})$ represents a two-dimensional subspace of $\mathcal{K}$. Each such subspace generates the space $S^1 \times H^3$ of branched immersions (equation (9.8)). □

<u>THEOREM 21.</u> *For each $p \geq 2$ there exists a real $4(p-1)$-dimensional family of minimal branched immersions of spheres punctured at $2p$ points with finite total curvature and embedded planar ends. For $2 \leq p \leq 7$, the moduli space of such immersions is exactly $4(p-1)$-dimensional.*

*Proof.* Let $P = \sum [a_i]$ be a divisor of degree $2p$ on $S^2$, and $S$ the unique spin structure on $S^2$. Let $\mathcal{H}$ and $\mathcal{K}$ be as in equation (14.11). Then pfaffian $\Omega = 0$ (see Appendix A) if and only if $\dim \mathcal{K} \geq 2$ if and only if there exists a surface with $2p$ ends at supp $P$. Counting real dimensions, the space of $2p$ ends is $4p$-dimensional; the Möbius transformations of $S^2$ reduce the dimension by 6, and the pfaffian condition on the ends reduce the dimension by another 2, so the space of ends which admit surfaces is $(4p - 8)$-dimensional. For each admissible choice of ends, by the above theorem there is a real $(4 \dim \mathcal{K} - 4)$-dimensional space of surfaces. Altogether, this totals $4p + 4 \dim \mathcal{K} - 12$, which is at least $4p - 4$ since $\dim \mathcal{K} \geq 2$.

In the case that $2 \leq p \leq 7$, by Lemma 19, $2 \leq \dim \mathcal{K} \leq [\sqrt{2p}] \leq [\sqrt{14}] = 3$, so $\dim \mathcal{K}$, being even, must be exactly 2. □

## 20  $\Omega$ on the Riemann sphere

For the examples in sections 21–23 we need to compute $\Omega$ on the Riemann sphere. Let $z$ be the standard conformal coordinate on $S^2 = \mathbb{C} \cup \{\infty\}$, and let $\varphi^2 = dz$ represent the unique spin structure on $S^2$. Let $P = [a_1] + \ldots + [a_{n-1}] + [\infty]$ with the $a_i \in \mathbb{C}$ distinct. We have $\mathcal{H} = 0$ since there are no holomorphic differentials on the sphere. A basis for $\mathcal{F}$ is

$$\{t_1, \ldots, t_{n-1}, t_n\} = \left\{ \frac{\varphi}{z - a_1}, \ldots, \frac{\varphi}{z - a_{n-1}}, \varphi \right\}.$$

These sections are in $\mathcal{F}$ since

$$(t_n) = -[\infty], \quad (t_i) = -[a_i],$$

and are independent because they have distinct poles, and so are a basis for $\mathcal{F}$ since $\dim \mathcal{F} = n$. By the local calculation (16.12) for $\Omega$,

$$\Omega(t_i, t_j) = \begin{cases} \dfrac{1}{a_j - a_i} & (1 \le i \le n-1;\ 1 \le j \le n-1;\ i \ne j), \\ -1 & (1 \le i \le n-1;\ j = n), \\ 1 & (i = n;\ 1 \le j \le n-1), \\ 0 & (i = j). \end{cases}$$

## 21 Genus zero surfaces with four embedded planar ends

The family of minimal genus zero surfaces with four embedded planar ends was computed first in [2]. A different computation is included here for completeness.

<u>THEOREM 22.</u>  *The space $\Sigma_4$ of complete minimal immersions of spheres punctured at four points into space with finite total curvature and embedded planar ends is $S^1 \times H^3$.*

*Proof.* Let $z$ be the standard conformal coordinate on $S^1 = \mathbb{C} \cup \{\infty\}$. By a Möbius transformation of the Riemann sphere $S^2$, the ends can be normalized so that two of the ends are 0 and $\infty$ and the product of the other two is 1. Naming the normalized ends

$$\{a_1 = a, a_2 = 1/a, 0, \infty\},$$

the matrix for $\Omega$ in the basis

$$\left\{\frac{1}{z - a_1}, \frac{1}{z - a_2}, \frac{1}{z}, 1\right\}$$

is

$$\Omega = \begin{pmatrix} 0 & \frac{1}{a_2 - a_1} & -\frac{1}{a_1} & -1 \\ \frac{1}{a_1 - a_2} & 0 & -\frac{1}{a_2} & -1 \\ \frac{1}{a_1} & \frac{1}{a_2} & 0 & -1 \\ 1 & 1 & 1 & 0 \end{pmatrix}$$

(see section 20). The pfaffian of $\Omega$ (see Appendix A) computes to a nonzero multiple of

$$(a^2 - \sqrt{3}a + 1)(a^2 + \sqrt{3}a + 1).$$

This pfaffian must be zero in order for $\ker \Omega = \mathcal{K}$ to be at least two-dimensional and hence to generate surfaces. Setting this pfaffian to zero yields interchangeable solutions for $a$, one of which is

$$a = (\sqrt{3} + i)/2.$$

With $\varphi^2 = dz$ as usual, a basis for $\mathcal{K}$ is

$$t_1 = \left(\frac{\sqrt{3}z - 1}{z(z^2 - \sqrt{3}z + 1)}\right)\varphi \quad \text{and} \quad t_2 = \left(\frac{z(z - \sqrt{3})}{z^2 - \sqrt{3}z + 1}\right)\varphi,$$

the family of immersions is then given by $X = \operatorname{Re} F$, where

$$F = \int (s_1^2 - s_2^2, i(s_1^2 + s_2^2), 2s_1 s_2)$$

and

$$\begin{pmatrix} s_1 \\ s_2 \end{pmatrix} = Q \begin{pmatrix} t_1 \\ t_2 \end{pmatrix},$$

where $Q \in \mathbb{C}^* \times \operatorname{SL}(2, \mathbb{C})$. The surfaces are identical (up to a rotation or dilation in space) when $Q \in \mathbb{R}^* \times \operatorname{SU}(2)$. Thus a parameter space for this family of surfaces is $S^1 \times H^3$ (see section 9). That these surfaces are immersed is shown in the next section. □

## 22 Genus zero surfaces with six embedded planar ends

Herein is computed the family of minimal genus zero surfaces with six embedded planar ends.

<u>THEOREM 23.</u>   *The space $\Sigma_6$ of complete minimal immersions of spheres punctured at six points into space with finite total curvature and embedded planar ends is $V \times (S^1 \times H^3)$, where $V$ is an algebraic subvariety of $(\mathbb{CP}^1)^3$ with codimension 1.*

*Proof.* On the sphere $S^2 = \mathbb{C} \cup \{\infty\}$ with standard conformal coordinate $z$, the ends can be normalized so that two of the ends are at $0$ and $\infty$, and the product of the remaining four ends is $1$. With this normalization, let the ends be $\{a_1, a_2, a_3, a_4, 0, \infty\}$. Set

$$\begin{aligned} \sigma_1 &= a_1 + a_2 + a_3 + a_4, \\ \sigma_2 &= -(a_1 a_2 + a_1 a_3 + a_1 a_4 + a_2 a_3 + a_2 a_4 + a_3 a_4), \\ \sigma_3 &= a_1 a_2 a_3 + a_1 a_2 a_4 + a_1 a_3 a_4 + a_2 a_3 a_4. \end{aligned}$$

The matrix for $\Omega$ in the basis
$$\left\{\frac{1}{z-a_1}, \frac{1}{z-a_2}, \frac{1}{z-a_3}, \frac{1}{z-a_4}, \frac{1}{z}, 1\right\}$$

is

$$\Omega = \begin{pmatrix} 0 & \frac{1}{a_2-a_1} & \frac{1}{a_3-a_1} & \frac{1}{a_4-a_1} & -\frac{1}{a_1} & -1 \\ \frac{1}{a_1-a_2} & 0 & \frac{1}{a_3-a_2} & \frac{1}{a_4-a_2} & -\frac{1}{a_2} & -1 \\ \frac{1}{a_1-a_3} & \frac{1}{a_2-a_3} & 0 & \frac{1}{a_3-a_4} & -\frac{1}{a_3} & -1 \\ \frac{1}{a_1-a_4} & \frac{1}{a_2-a_4} & \frac{1}{a_3-a_4} & 0 & -\frac{1}{a_4} & -1 \\ \frac{1}{a_1} & \frac{1}{a_2} & \frac{1}{a_3} & \frac{1}{a_4} & 0 & -1 \\ 1 & 1 & 1 & 1 & 1 & 0 \end{pmatrix}$$

(see section 9). The pfaffian of $\Omega$ (see Appendix A) is

(22.14) $$\text{pfaffian}\,\Omega = \tau_1\tau_3 + \sigma_1\sigma_3 - 20,$$

where
$$\tau_1 = \sigma_1^2 + 3\sigma_2 \quad \text{and} \quad \tau_3 = \sigma_3^2 + 3\sigma_2.$$

The condition that the pfaffian be 0 defines the algebraic subvariety

$$V = \{(\sigma_1, \sigma_2, \sigma_3) \subset (\mathbb{CP}^1)^3 \mid \text{pfaffian}\,\Omega = 0\}$$

of codimension 1. Assuming that the pfaffian is zero, and letting $\varphi^2 = dz$, a basis for the kernel of $\Omega$ is $\{t_1, t_2\}$, where

$$t_1 = \left(\frac{b_3 z^3 + b_2 z^2 + b_1 z + b_0}{z(z^4 - \sigma_1 z^3 - \sigma_2 z^2 - \sigma_3 z + 1)}\right)\varphi,$$
$$t_2 = \left(\frac{z(c_3 z^3 + c_2 z^2 + c_1 z + c_0)}{z^4 - \sigma_1 z^3 - \sigma_2 z^2 - \sigma_3 z + 1}\right)\varphi,$$

and

$$\begin{aligned}
b_0 &= \sigma_2, & c_0 &= \sigma_3\tau_1 + 5\sigma_1, \\
b_1 &= -\sigma_2\sigma_3, & c_1 &= \sigma_2\tau_1 - 2\sigma_1\sigma_3 - 10, \\
b_2 &= \sigma_2\tau_3 - 2\sigma_1\sigma_3 - 10, & c_2 &= -\sigma_1\sigma_2, \\
b_3 &= \sigma_1\tau_3 + 5\sigma_3, & c_3 &= \sigma_2
\end{aligned}$$

(the special case $\sigma_2 = 0$ for which the above sections are linearly dependent is ignored here). The family of immersions is then given by

$$X = \operatorname{Re} \int (s_1^2 - s_2^2, i(s_1^2 + s_2^2), 2s_1 s_2)$$

where

$$\begin{pmatrix} s_1 \\ s_2 \end{pmatrix} = Q \begin{pmatrix} t_1 \\ t_2 \end{pmatrix},$$

and $Q \in \mathbb{C}^* \times \operatorname{SL}(2, \mathbb{C}) / \mathbb{R}^* \times \operatorname{SU}(2) \cong S^1 \times H^3$ as in the previous section. □

That the four- and six-ended families are immersed follows from the lemma below, which in turn follows directly from the definitions of the spaces in equation (14.10).

<u>LEMMA 24.</u> *On the sphere with its unique spin structure $S$, let $P_1 = \sum [p_i]$ as in equation (14.10), and $P_2 = P_1 + [a]$, $(a \notin \operatorname{supp}(P_1))$. Let $\mathcal{F}_i = \mathcal{F}_{S^1, P_i, S}$ and $\mathcal{K}_i = \mathcal{K}_{S^1, S, P_i}$ $(i = 1, 2)$ as in equation (14.11). Then $\mathcal{K}_2 \cap \mathcal{F}_1 = \{s \in \mathcal{K}_1 \mid s(a) = 0\}$.*

Now, to complete the proof that the above examples are immersed, let $P_1$ be the divisor of ends of even degree $n < 9$, and let $(s_1, s_2)$ be the spinor representation of a minimal branched immersion. Supposing this surface is not immersed, let $a$ be a branch point of the surface, and set $P_2 = P_1 + [a]$. Then $s_1$ and $s_2$ are independent sections in $\mathcal{K}_1$ and $s_1(a) = 0$, $s_2(a) = 0$, so by the lemma, $s_1, s_2 \in \mathcal{K}_2$. Applying Lemma 19 iii), we have that

$$2 \leq \dim \mathcal{K}_2 \leq [\sqrt{n}] \leq 2,$$

so $\dim \mathcal{K}_2 = 2$. This contradicts the fact that $n + 1 - \dim \mathcal{K}_2$ is even (Lemma 19 i)).

## 23 Projective planes with three embedded planar ends

It was shown in [14] that any minimal immersion of a punctured real projective plane with embedded ends has only planar ends, and has at least three of them. Hence those which are the subject of the following theorem are the examples of minimal projective planes with the fewest number of embedded ends. One method for determining the moduli space of finite total curvature minimally immersed projective planes punctured at three points was given in [3]. Here we provide another description of this moduli space using the spinor representation. Note that all these surfaces compactify to give surfaces minimizing $W = \int H^2 dA$ among all immersed real projective planes [13], with minimum energy $W = 12\pi$.

THEOREM 25. *Let $\Pi_3$ be the moduli space of complete minimal immersions of real projective planes punctured at three points with finite total curvature and embedded planar ends modulo Euclidean similarities. Then*
  (i) *$\Pi_3$ is homeomorphic to a closed disk with one point $M_0$ removed from the boundary;*
 (ii) *the point $M_0$ represents the Möbius strip with total curvature $-6\pi$ in the sense that if $\gamma : \mathbb{R}^+ \longrightarrow \Pi_3$ is a curve with $\lim_{t\to\infty} \gamma(t) = M_0$, then there is a one-parameter family of immersions $X_t$ parametrizing the surfaces $\gamma(t)$ such that as $t \to \infty$, $X_t$ converges uniformly on compact sets to a parametrization of the Möbius strip;*
(iii) *the surfaces with non-trivial symmetry groups are represented by the boundary of the disk, which represents a one-parameter family of surfaces which have a line of reflective symmetry; among these, the only surfaces with larger symmetry groups (other than $M_0$) are two surfaces which have, respectively, the symmetry groups $\mathbb{Z}_2 \times \mathbb{Z}_2$, and $D_3$, the dihedral group of order 6.*

*Proof of (i):* The two-sheeted covering of the projective plane is the Riemann sphere $S^2 = \mathbb{C} \cup \{\infty\}$, with order-two orientation-reversing deck transformation $I(z) = -1/\overline{z}$. By a motion in PSU(2) the six preimages on the sphere of three points in the projective plane can be normalized as in section 22 to be

$$\{a_1, I(a_1), a_2, I(a_2), 0, \infty\}$$

with the product of the first four equal to 1. With this choice, following the notation of section 22, we have
$$\sigma_2 \in \mathbb{R}; \quad \sigma_3 = -\overline{\sigma}_1; \quad \tau_3 = \overline{\tau}_1.$$

For each choice of ends satisfying equation (22.14), up to dilations and isometries of space there is a unique minimal immersion of the projective plane, whose spinor representation is given by $\sqrt{i}(t_1, t_2)$, with $t_1, t_2$ as in section 22. For if $\sqrt{i}(\hat{t}_1, \hat{t}_2)$ is the spinor representation of another immersion with the same ends, then a motion in $\mathbb{C}^* \times \text{PSL}(2,\mathbb{C})$ can make $\hat{t}_1 = t_1$, and the compatibility condition in Theorem 11 forces $\hat{t}_2 = \pm t_1$. Hence the moduli space $\Pi_3$ can be parametrized as a quotient space of
$$\Gamma = \{(\sigma_1, \sigma_2) \in \mathbb{C} \times \mathbb{R} \mid \tau_1 \tau_3 + \sigma_1 \sigma_3 - 20 = 0, \ \sigma_3 = -\overline{\sigma}_1\},$$
where $\sigma_1, \sigma_2, \sigma_3$ are the symmetric polynomials of the ends defined in section 22. The desired moduli space is a quotient space of $\Gamma$, since permutations of the ends give rise to the same surface.

Since the parameters $\sigma_1$ and $\sigma_2$ depend on the particular normalization of the ends made in section 22, new parameters should be chosen, namely the three direction cosines of the angles between the ends $0$, $a_1$ and $a_2$, viewed as vectors in $S^2 \subset \mathbb{R}^3$. To convert the equation (22.14) to these new parameters let $\phi : \mathbb{C} \longrightarrow S^2 \subset \mathbb{R}^3$ be inverse stereographic projection defined by

$$\phi(a) = \left( \frac{2\operatorname{Re} a}{|a|^2 + 1}, \frac{2\operatorname{Im} a}{|a|^2 + 1}, \frac{|a|^2 - 1}{|a|^2 + 1} \right).$$

With the usual inner product $\langle \, , \, \rangle$ in $\mathbb{R}^3$, the direction cosines are

$$c_1 = \langle \phi(0), \phi(a_1) \rangle = \frac{1 - |a_1|^2}{1 + |a_1|^2},$$

$$c_2 = \langle \phi(0), \phi(a_2) \rangle = \frac{1 - |a_2|^2}{1 + |a_2|^2},$$

$$c_3 = \langle \phi(a_1), \phi(a_2) \rangle = \frac{(1 - |a_1|^2)(1 - |a_2|^2)}{(1 + |a_1|^2)(1 + |a_2|^2) + 4\operatorname{Re} a_1 \bar{a}_2}.$$

The above three equations may be written

$$|a_1|^2 = \frac{1 - c_1}{1 + c_1},$$

$$|a_2|^2 = \frac{1 - c_2}{1 + c_2},$$

$$\operatorname{Re} a_1 \bar{a}_2 = \frac{c_3 - c_1 c_2}{(1 + c_1)(1 + c_2)}.$$

Using the normalization of the ends above, and writing $a_1 = \gamma r_1$, $a_2 = \bar{\beta} r_2$ ($\gamma \in S^1 \subset \mathbb{C}$; $r_1, r_2 \in \mathbb{C}$) yields

$$r_1 = \frac{\sqrt{1 - c_1}}{\sqrt{1 + c_1}},$$

$$r_2 = \frac{\sqrt{1 - c_2}}{\sqrt{1 + c_2}},$$

$$\gamma^2 = \frac{c_3 - c_1 c_2 + ix}{\sqrt{1 - c_1^2}\sqrt{1 - c_2^2}},$$

where

$$x^2 = 1 - c_1^2 - c_2^2 - c_3^2 + 2c_1 c_2 c_3.$$

To convert the determinant of equation (22.14) from the variables $\sigma_1$, $\sigma_2$, $\sigma_3$ to $c_1$, $c_2$, $c_3$, compute

$$\sigma_1 = \frac{2\gamma(-c_1 + 2c_1 c_2^2 - c_2 c_3 + ic_2 x)}{\sqrt{1 - c_1^2}(1 - c_2^2)},$$

$$\sigma_2 = \frac{2(c_3 - 3c_1 c_2)}{(1 - c_1^2)(1 - c_2^2)},$$

$$\sigma_3 = -\overline{\sigma}_1,$$

and the determinant becomes, up to a non-zero multiple,

$$(c_1^2 + 3)(c_2^2 + 3)(c_3^2 + 3) - 32(c_1 c_2 c_3 + 1).$$

The surface

$$\Gamma = \left\{(c_1, c_2, c_3) \in \mathbb{R}^3 \mid (c_1^2 + 3)(c_2^2 + 3)(c_3^2 + 3) - 32(c_1 c_2 c_3 + 1) = 0\right\}$$

in the cube

$$C = \left\{(x, y, z) \in \mathbb{R}^3 \mid -1 < x, y, z < 1\right\}$$

is a tetrahedron-like object but with smoothed edges and (omitted) vertices at $(1, 1, 1)$, $(1, -1, -1)$, $(-1, 1, -1)$, and $(-1, -1, 1)$.

The moduli space $\Pi_3$ is a quotient of $\Gamma$ which arises from permutations of the ends. A choice $c = (c_1, c_2, c_3)$ determines a set of six ends on the double-covering sphere. The group of rotations of the cube is the order-24 permutation group $S_4$ generated by two kinds of elements:

- permuting the three numbers $(c_1, c_2, c_3)$,
- negating any two of the three numbers $(c_1, c_2, c_3)$.

Action under this group determines the same six ends. Hence $\Pi_3 = \Gamma/S_4$ is a representation of the moduli space of minimal projective planes with three embedded planar ends.

Draw the two diagonals on each face of the cube $C$ dividing each face into four triangles. Consider the the 24 tetrahedra whose bases are these triangles, and whose common vertex is the origin. Each of these tetrahedra is a fundamental domain under the action of $S_4$ on the cube. This can also be seen by noting that any $(c_1', c_2', c_3')$ in the cube $C$ has in its orbit under $S_4$ a point $(c_1, c_2, c_3)$ satisfying $c_1 \geq c_2 \geq |c_3| \geq 0$. Let

$$T = \{(c_1, c_2, c_3) \in C \mid c_1 \geq c_2 \geq |c_3| \geq 0\}$$

be one of these tetrahedra. Then $D = T \cap \Gamma$ is a fundamental domain in $\Gamma$ for the group $S_4$, with boundary

$$\partial D = \partial T \cap \Gamma = (\{c_1 = c_2\} \cup \{c_2 = c_3\} \cup \{c_2 = -c_3\}) \cap (T \cap \Gamma).$$

$D$ can be shown to be topologically a closed disk with the point corresponding to the corner $(1,1,1)$ of the cube removed.

*Proof of (ii):* The minimal Möbius strip with total curvature $-6\pi$, found in [18], has spinor representation

$$G(w)\sqrt{dw} = \sqrt{i}(-(w+1)/w^2, w-1)\sqrt{dw}$$

Let $(\sigma_1(s), \sigma_2(s)) : \mathbb{R}^+ \longrightarrow \Gamma$ be a proper curve. It follows from the reality of $\sigma_2$ that

$$\lim_{s\to\infty} \frac{1}{\sigma_1(s)} = \lim_{s\to\infty} \frac{1}{\sigma_2(s)} = \lim_{s\to\infty} \frac{\sigma_1(s)}{\sigma_2(s)} = 0,$$

and by a permutation of the ends we can assume

$$\lim_{s\to\infty} \frac{\overline{\sigma_1(s)}}{\sigma_1(s)} = 1.$$

Further,

$$\lim_{s\to\infty} \left|\frac{\tau_1(s)}{\sigma_1(s)}\right| = 1,$$

since

$$\left|\frac{\tau_1}{\sigma_1}\right|^2 = -\frac{\tau_1\tau_3}{\sigma_1\sigma_3} = 1 - \frac{20}{|\sigma_1^2|}.$$

Now choose a function $\alpha : \mathbb{R}^+ \longrightarrow S^1 \subset \mathbb{C}$ such that

$$\lim_{s\to\infty}\left(\frac{\tau_1(s)}{\sigma_1(s)} - \alpha(s)\right) = 0,$$

and so

$$\lim_{s\to\infty}\left(\frac{\tau_3(s)}{\sigma_1(s)} - \overline{\alpha(s)}\right) = 0.$$

Let $X$ be defined by

$$X(z)\sqrt{dz} = \frac{\sqrt{i}}{\sigma_1}(t_1, t_2),$$

where $t_1$, $t_2$ are as in section 22. A careful reparametrization and rotation of the surface generated by $X(z)\sqrt{dz}$ converges uniformly in compact sets to the Möbius strip given above: Let $z = \alpha w$, and

$$A_\alpha = \begin{pmatrix} \alpha^{3/2} & 0 \\ 0 & \alpha^{-3/2} \end{pmatrix}.$$

Then

$$A_\alpha X(z)\sqrt{dz} = A_\alpha \sqrt{\alpha} X(\alpha w)\sqrt{dw}$$

is the appropriate reparametrization and rotation. This amounts to showing

$$\lim_{s\to\infty} A_{\alpha(s)}\sqrt{\alpha(s)}X(\alpha(s)w) = G(w)$$

uniformly in compact sets not containing the ends, which follows by a calculation using the limits above.

*Proof of (iii):* To find the surfaces in $\Pi_3$ which have non-trivial symmetry groups as surfaces in space, let $G = \mathbb{Z}_2 \times \text{PSU}(2) \cong O(3)$ be the group of conformal and anticonformal diffeomorphisms of $\mathbb{C} \cup \{\infty\} = S^2$ with the property that any $\xi \in G$ commutes with $I$. Via stereographic projection, $G$ can be thought of as the isometry group of $S^2 \subset \mathbb{R}^3$, so $\xi \in G$ satisfies $\langle a, b \rangle = \langle \xi a, \xi b \rangle$. The group of symmetries of the minimal surface in space induces a subgroup $H \subset G$ acting on the domain $S^2$. Moreover, the subgroup $H \subset G$ which permutes the ends is isomorphic to the subgroup $K \subseteq S_4$ which fixes the point $(c_1, c_2, c_3)$ representing the ends, since $\xi \in H$ preserves the inner product defining the cosines $c_1$, $c_2$, $c_3$.

The point of all this is that the symmetry group of a surface represented by $(c_1, c_2, c_3) \in \Pi_3$ can be determined by finding the subgroup of $S_4$ which fixes $(c_1, c_2, c_3)$. Using this method, the surfaces other than the Möbius strip at $(1, 1, 1)$ are

- elements of $\partial D$, each with a line of reflective symmetry,
- $(\sqrt{5}/3, 0, 0) \in \partial D$ with symmetry group $\mathbb{Z}_2 \times \mathbb{Z}_2$,
- $(c, c, -c) \in \partial D$ with symmetry group $S_3 = D_3$

The last (and most symmetric) of these is a surface described in [14]. □

## 24 Genus one

The remaining sections concern minimal immersions in the regular homotopy classes of tori and Klein bottles with embedded planar ends. In sections 25 and 26, the skew-symmetric form $\Omega$ is computed for the twisted and the untwisted tori. This computation is then used to show the nonexistence and existence of various examples. In section 27 it is shown that no such tori exist with three ends, and in section 28, is found a real two-dimensional family of immersions with four ends exists on each conformal type of torus. After some general results about Klein bottles in section 29, a minimal Klein bottle with embedded planar ends is constructed in section 30.

## 25 $\Omega$ on the twisted torus

For the non-example in section 27, and for the example in section 28, it is necessary to compute a basis for $\mathcal{F}$ for the twisted torus (see section 8), and the matrix for $\Omega$ in this basis. On the torus $\mathbb{C}/\{2\omega_1, 2\omega_3\}$ with the standard coordinate $u$, let $S$ be the spin structure corresponding to the twisted torus, that is, represented by the holomorphic differential $\varphi_0^2 = du$. Let $P = [a_1] + \ldots + [a_n]$ and set $\omega_2 = \omega_1 + \omega_3$ throughout the remainder of Part III.

To show that $\mathcal{H} = \{c\varphi_0 \mid c \in \mathbb{C}\}$, let $t \in \mathcal{H}$. Then $0 \leq (t) = (t/\varphi_0) + (\varphi_0) = (t/\varphi_0)$. Hence $t/\varphi_0$ is a holomorphic function on the torus, so it is constant.

A basis for $\mathcal{F}$ is $\{t_0, t_1, \ldots, t_{n-1}\}$, where

$$t_0 = \varphi_0,$$
$$t_i = (\zeta(u - a_i) - \zeta(u) + \zeta(a_i))\,\varphi_0,$$
$$= \frac{1}{2}\left(\frac{\wp'(u) + \wp'(a_i)}{\wp(u) - \wp(a_i)}\right)\varphi_0$$

(see equation (B.17)). These are in $\mathcal{F}$ because

$$(t_0) = 0 \geq -P,$$
$$(t_i) = [x_i] + [y_i] - [a_i] - [0] \geq -P$$

where $x_i$ and $y_i$ are the zeros of $\wp'(u) + \wp'(a_i)$ other than $-a_i$. These sections are independent because they have distinct poles, and hence span $\mathcal{F}$ since $\dim \mathcal{F} = n$. To compute $\Omega$ in this basis, first compute the expansions of $t_i$ at $a_0, \ldots, a_{n-1}$ (assume $i$, $j \neq 0$):

$$t_i = (-u^{-1} + \mathrm{o}(u))\varphi_0,$$
$$t_i = ((t_i/\varphi_0)(a_j) + \mathrm{o}(u))\varphi_0 \quad (i \neq j),$$
$$t_i = (u - a_i)^{-1}\varphi_0.$$

Using equation (16.12), we have

$$\Omega(t_i, t_j) = \begin{cases} \left.\dfrac{t_i}{\varphi_0}\right|_{a_j} & (i \neq 0;\ j \neq 0;\ i \neq j), \\ 0 & (\text{otherwise}). \end{cases}$$

## 26 $\Omega$ on the untwisted tori

As above, it is also necessary to exhibit a basis for $\mathcal{F}$ on the untwisted tori (see section 8), as well as the matrix for $\Omega$ in this basis. On the torus $\mathbb{C}/\{2\omega_1, 2\omega_3\}$ with the standard conformal coordinate $u$, fix $r \in \{1, 2, 3\}$ and choose the spin structure on the untwisted torus, represented by

$$\varphi_r^2 = \frac{du}{\wp_r(u)},$$

where $\wp_r(u) = \wp(u) - \wp(\omega_r)$. Let $P = \sum[a_i]$ with the $a_i \in T \setminus \{0, \omega_r\}$ distinct.

For this choice of spin structure, $\mathcal{H} = 0$. To show this, first note first that $(\varphi_r) = [0] - [\omega_r]$. If $t \in \mathcal{H}$, then

$$0 \leq (t) = (t/\varphi_r) + (\varphi_r) = (t/\varphi_r) + [0] - [\omega_r].$$

It follows that $(t/\varphi_r) \geq [\omega_r] - [0]$. But since $t/\varphi_r$ is a function, the degree of its divisor is 0. Hence $(t/\varphi_r) = [\omega_r] - [0]$. But this is impossible by Abel's theorem on the torus: for an elliptic function $f$, if $(f) = \sum n_i[p_i]$ (as a formal sum) then $\sum n_i p_i = 0$ (as a sum in $\mathbb{C}$).

A basis for $\mathcal{F}$ is $\{t_1, \ldots, t_n\}$, where

$$t_i(u) = (\zeta(u - a_i) - \zeta(u) - \zeta(\omega_r - a_i) + \zeta(\omega_r))\varphi_r$$
$$= \frac{1}{2}\left(\frac{\wp_r(u)\wp'_r(a_i) + \wp'_r(u)\wp_r(a_i)}{\wp_r(a_i)(\wp_r(u) - \wp_r(a_i))}\right)\varphi_r$$

(see equation (B.17)). These are in $\mathcal{F}$ because $(\varphi_r) = [0] - [\omega_r]$, so $(t_i) = [a_i - \omega_r] - [a_i] \geq -P$, and are independent because their poles are distinct, so they span $\mathcal{F}$ since $\dim \mathcal{F} = n$. The expansions of $t_i$ at $a_1, \ldots, a_n$ are

$$t_i = ((t_i/\varphi_r)(a_j) + o(u - a_j))\varphi_r \quad (i \neq j),$$
$$t_i = ((u - a_i)^{-1} + o(u - a_i))\varphi_r.$$

Using the local expression (16.12) for $\Omega$, we have

$$\Omega(t_i, t_j) = \begin{cases} \left.\dfrac{t_i}{\varphi_r}\right|_{a_j} & (i \neq j), \\ 0 & (i = i). \end{cases}$$

A particularly simple situation arises when the ends come in pairs $a$ and $-a$. Assume $n = 2m$ and $a_{m+i} = -a_i$ ($i = 1, \ldots, m$). In this case, a simpler basis is $\{\hat{t}_1, \ldots, \hat{t}_m, \hat{t}_{m+1}, \ldots, \hat{t}_{2m}\}$, where for $1 \leq i \leq m$,

$$\hat{t}_i = \frac{\wp_r(a_i)}{\wp'_r(a_i)}(t_i - t_{m+i})\varphi_r = \left(\frac{\wp_r(u)}{\wp_r(u) - \wp_r(a_i)}\right)\varphi_r,$$

$$\hat{t}_{m+i} = \qquad (t_i + t_{m+i})\varphi_r = \left(\frac{\wp'_r(u)}{\wp_r(u) - \wp_r(a_i)}\right)\varphi_r.$$

In this basis, the matrix for $\Omega$ becomes

$$\left(\begin{array}{c|c} 0 & W \\ \hline -W^t & 0 \end{array}\right),$$

where $W$ is given by

$$W_{ij} = \begin{cases} \dfrac{4}{\wp_r(a_i) - \wp_r(a_j)} & (i < j), \\ \dfrac{4}{\wp_r(a_j) - \wp_r(a_i)} & (i > j), \\ \dfrac{\wp_r(a_i)^2 - c_p c_q}{\wp_r(a_i)(\wp_r(a_i) - c_p)(\wp_r(a_i) - c_q)} & (i = j) \end{cases}$$

and $c_p = e_p - e_r$, $c_q = e_q - e_r$, $\{p, q, r\} = \{1, 2, 3\}$. Note that the entries of $W$ are entirely free of $\wp_r'$.

A useful property of the basis above is as follows: let $L : M \longrightarrow M$ be defined as $L(u) = -u$; then for $i \leq m$ and $j \geq m+1$,

$$L^*(\hat{t}_i \hat{t}_j) = \hat{t}_i \hat{t}_j,$$

so

$$\int_{\gamma_k} \hat{t}_i \hat{t}_j = \int_{\gamma_k} L^*(\hat{t}_i \hat{t}_j) = \int_{L(\gamma_k)} \hat{t}_i \hat{t}_j = -\int_{\gamma_k} \hat{t}_i \hat{t}_j.$$

and so

$$\int_{\gamma_k} \hat{t}_i \hat{t}_j = 0 \quad (i \leq m; \ j \geq m+1; \ k = 1, 3).$$

# 27  Non-existence of tori with three planar ends

An outline of the proof of the non-existence of three-ended tori, twisted or untwisted, is given.

THEOREM 26.  *There does not exist a complete minimal branched immersion of a torus into space with finite total curvature and three embedded planar ends.*

*Sketch of proof:* The proof is divided into two cases: for the twisted torus there exist immersions with periods, but the periods cannot be made purely imaginary; for the untwisted torus, there are not even periodic examples.

First consider the more difficult case of of the twisted torus. With everything as in section 25, let $\{0, a_1, a_2\}$ be the set of ends, and let $p_i = \wp(a_i)$, $p_i' = \wp'(a_i)$. The condition $\dim \mathcal{K} \geq 2$ puts the following condition on the placement of the ends:

$$g_2 = 4(p_1^2 + p_1 p_2 + p_2^2),$$

where $g_2$ is the constant in the differential equation $(\wp')^2 = 4\wp^3 - g_2\wp - g_3$. To see this, first note that $\ker \Omega = \mathcal{K} \oplus \mathcal{H}$ and $\dim \mathcal{H} = 1$. Hence if $\dim \mathcal{K} = 2$ then $\Omega \equiv 0$.

Assume first that $a_1 + a_2 \neq 0$. Then $p_1 - p_2 \neq 0$, and the entries of $\Omega$ indicate that $p'_1 + p'_2 = 0$ Hence
$$(p'_1)^2 = 4p_1^3 - g_2 p_1 - g_3$$
and
$$(p'_2)^2 = 4p_2^3 - g_2 p_2 - g_3$$
are equal, and the desired condition follows. The condition also obtains in the case that $a_1 + a_2 = 0$; this can be shown as a limiting case of the above.

Changing basis now to simplify the period equations, let
$$\hat{t}_1 = t_1 + \varepsilon t_2,$$
$$\hat{t}_2 = t_1 + \varepsilon^2 t_2,$$
where $\varepsilon = (-1 + \sqrt{3})/2$. With $\gamma_1$, $\gamma_3$ the closed curves parallel to $\omega_1$, $\omega_3$ respectively (as in Theorem 28), the integrals relevant to the periods are (for $k = 1, 3$)
$$\int_{\gamma_k} \hat{t}_1^2 = -6q_1 \omega_k,$$
$$\int_{\gamma_k} \hat{t}_1 \hat{t}_2 = -6\eta_k,$$
$$\int_{\gamma_k} \hat{t}_2^2 = -6q_2 \omega_k,$$
where
$$q_1 = -((\varepsilon - \varepsilon^2)p_1 + (\varepsilon - 1)p_2)/3,$$
$$q_2 = -((\varepsilon^2 - \varepsilon)p_1 + (\varepsilon^2 - 1)p_2)/3,$$
$$q_1 q_2 = (p_1^2 + p_1 p_2 + p_2^2)/3 = g_2/12.$$

A choice of a pair of independent sections from $\mathcal{K}$ can be normalized by the action of $\mathbb{R}^* \times \mathrm{SU}(2)$ to be
$$s_1 = z_1 \hat{t}_1 + \hat{t}_2,$$
$$s_2 = z_2 \hat{t}_1,$$
with $z_1, z_2 \in \mathbb{C}$. Then the period equations (10.9) can be written
$$\begin{pmatrix} 2z_1 \\ z_1^2 q_1 + q_2 \end{pmatrix} - B \begin{pmatrix} 0 \\ q_1 \bar{z}_2^2 \end{pmatrix} = 0,$$
$$\begin{pmatrix} z_2 \\ q_1 z_1 z_2 \end{pmatrix} + B \begin{pmatrix} \bar{z}_2 \\ q_1 \bar{z}_1 \bar{z}_2^2 \end{pmatrix} = 0,$$

where
$$B = A^{-1}\overline{A} = \begin{pmatrix} a & b \\ c & d \end{pmatrix}; \quad A = \begin{pmatrix} \eta_1 & \omega_1 \\ \eta_3 & \omega_3 \end{pmatrix}.$$

Changing from the variables $(z_1, z_2)$ to $(w, z_2)$, this system is equivalent to the system

$$w^2 + b^2 q_1 q_2 - d^2 = 0,$$
$$2w + 2d - b^2 q_1 \overline{q}_1 \overline{z}_2^2 = 0,$$
$$wz_2 + \overline{z}_2 = 0.$$

From these it follows that

$$w\overline{w} - 1 = 0,$$
$$aw^2 - \overline{a} = 0,$$
$$-\overline{a} - ab^2 q_1 q_2 + ad^2 = 0.$$

This last condition, depending only on the conformal type of the torus and not on $w$, $z_1$, and $z_2$, is a degeneracy condition for the period equations. It also follows, by an examination of the sign of $a(w - \overline{a}) \in \mathbb{R}$, that

$$|a| > 1.$$

A delicate argument, which is omitted here, using the expansions [16]

$$g_2 = \frac{\pi^4}{12\omega_1^4}\left(1 + 240 \sum_{n=1}^{\infty} \sigma_3(n) q^n\right),$$
$$\eta_1 = \frac{\pi^2}{12\omega_1}\left(1 - 24 \sum_{n=1}^{\infty} \sigma_1(n) q^n\right),$$

where

$$\sigma_k(n) = \sum_{d|n} d^k; \quad q = e^{2i\pi\tau}; \quad \tau = \omega_3/\omega_1$$

shows that the degeneracy condition is not satisfied under the constraint $|a| > 1$ over the whole moduli space of Riemann tori. Hence no examples with three ends can be found in the case of the twisted tori.

The case of the untwisted tori is much easier. Fix $r \in \{1, 2, 3\}$ and let $\varphi_r$ be as in section 26. Let $\{a_1, a_2, a_3\}$ be the ends, translated so that they avoid $\{0, \omega_r\}$, and let $\{t_1, t_2, t_3\}$ be the basis for $\mathcal{F}$ given in the same section. The condition that $\dim \mathcal{K} = \dim \ker \Omega \leq 2$ forces $\Omega$ to be zero. This means, for example, that $t_1/\varphi_r$ have zeros at $a_2$ and $a_3$. But the zeros of $t_1/\varphi_r$ are $\omega_r$ and $a_1 - \omega_r$, so one of $a_2, a_3$ has to be $\omega_r$, contrary to the assumption. □

## 28  Minimal tori with four embedded planar ends

Here the existence of families of four-ended tori is established.

<u>Theorem 27.</u>  *For each conformal type of torus there exists a real two-dimensional family of complete minimal immersions of the torus punctured at four points into space with finite total curvature and embedded planar ends. Each of the tori is twisted.*

*Proof.*  To exhibit the family, it is first necessary to determine the placement of the four ends. The ends in fact must be, up to a translation, at the four half-lattice points. To show this, on the torus $\mathbb{C}/\{2\omega_1, 2\omega_3\}$, assume the four ends are $\{0, a_1, a_2, a_3\}$, where $a_1, a_2, a_3$ are distinct points in the torus to be determined. With $\varphi_0^2 = du$, the matrix for $\Omega$ in the basis $\{1, t_1, t_2, t_3\} = \{\varphi, f_1\varphi_0, f_2\varphi_0, f_3\varphi_0\}$ of section 25 is

$$\Omega = \begin{pmatrix} 0 & 0 & 0 & 0 \\ 0 & 0 & f_1(a_2) & f_1(a_3) \\ 0 & f_2(a_1) & 0 & f_2(a_3) \\ 0 & f_3(a_1) & f_3(a_2) & 0 \end{pmatrix}.$$

If $\ker \Omega = \mathcal{H} \oplus \mathcal{K}$ is two-dimensional, then $\dim \mathcal{K} = 1$, since $\dim \mathcal{H} = 1$, so $\mathcal{K}$ is not big enough to generate a minimal surface. Hence to produce surfaces, $\operatorname{rank} \Omega$, being even, must be zero. In this case, all the entries of the above matrix are zero; a look at $t_i$ shows that $\wp'(a_i) + \wp'(a_j) = 0$ for all $i \neq j$. It follows that $\wp'(a_1) = \wp'(a_2) = \wp'(a_3) = 0$, so $\{a_1, a_2, a_3\} = \{\omega_1, \omega_2, \omega_3\}$.

With the ends fixed at $\{0, \omega_1, \omega_2, \omega_3\}$, $\mathcal{F} = \ker \Omega = \mathcal{H} \oplus \mathcal{K}$, so $\{t_1, t_2, t_3\}$ is a basis for $\mathcal{K}$. The simple zeros and poles of $t_1$, $t_2$, and $t_3$ are illustrated below.

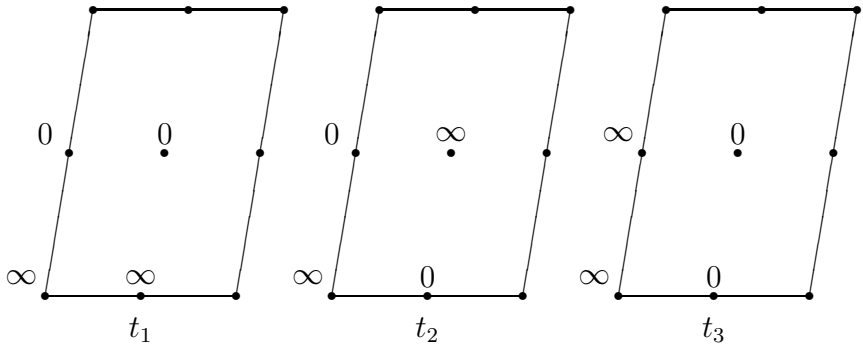

Figure 6: Zeros and poles of $t_1$, $t_2$, and $t_3$

To solve the period problem outlined in section 10 it is convenient to choose a new basis $\{\hat{t}_1, \hat{t}_2, \hat{t}_3\}$ for $\mathcal{K}$ which "diagonalizes" the period equations. Let

$$\begin{pmatrix} \hat{t}_1 \\ \hat{t}_2 \\ \hat{t}_3 \end{pmatrix} = \begin{pmatrix} 1 & -1 & -1 \\ -1 & 1 & -1 \\ -1 & -1 & 1 \end{pmatrix} \begin{pmatrix} t_1 \\ t_2 \\ t_3 \end{pmatrix},$$

or

$$\hat{t}_1(u) = (\zeta(u) + \zeta(u - \omega_1) - \zeta(u - \omega_2) - \zeta(u - \omega_3) + 2\zeta(\omega_1))\varphi_0,$$
$$\hat{t}_2(u) = (\zeta(u) - \zeta(u - \omega_1) + \zeta(u - \omega_2) - \zeta(u - \omega_3) + 2\zeta(\omega_2))\varphi_0,$$
$$\hat{t}_3(u) = (\zeta(u) - \zeta(u - \omega_1) - \zeta(u - \omega_2) + \zeta(u - \omega_3) + 2\zeta(\omega_3))\varphi_0.$$

The simple zeros and poles of $\hat{t}_1$, $\hat{t}_2$, and $\hat{t}_3$ are illustrated below. To compute the

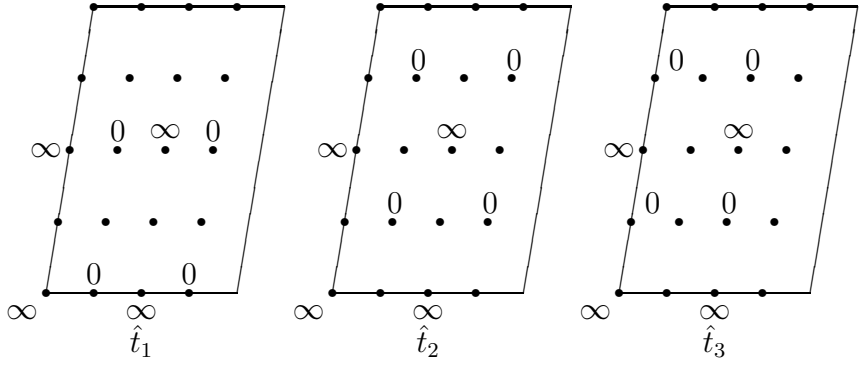

Figure 7: Zeros and poles of $\hat{t}_1$, $\hat{t}_2$, and $\hat{t}_3$

periods, use equation (B.18) to write

$$\hat{t}_i^2(u) = (\wp(u) + \wp(u - \omega_1) + \wp(u - \omega_2) + \wp(u - \omega_3) - 4\wp(\omega_i))\,du,$$
$$(\hat{t}_1\hat{t}_2)(u) = (\wp(u) - \wp(u - \omega_1) - \wp(u - \omega_2) + \wp(u - \omega_3))\,du,$$
$$(\hat{t}_1\hat{t}_3)(u) = (\wp(u) - \wp(u - \omega_1) + \wp(u - \omega_2) - \wp(u - \omega_3))\,du,$$
$$(\hat{t}_2\hat{t}_3)(u) = (\wp(u) + \wp(u - \omega_1) - \wp(u - \omega_2) - \wp(u - \omega_3))\,du.$$

With $\gamma_1$, $\gamma_3$ the closed curves on the torus respectively parallel to $\omega_1$, $\omega_3$, the periods are

$$P_k^{ij} = \int_{\gamma_k} \hat{t}_i \hat{t}_j du = \begin{cases} -8(\eta_k + \omega_k e_i) & \text{if } i = j \\ 0 & \text{if } i \neq j \end{cases} \quad (k = 1, 3),$$

where $e_i = \wp(\omega_i)$ and $\eta_k = \zeta(\omega_k)$ (see appendix B). In general, with

$$t_1 = x_1\hat{t}_1 + x_2\hat{t}_2 + x_3\hat{t}_3$$
$$t_2 = y_1\hat{t}_1 + y_2\hat{t}_2 + y_3\hat{t}_3$$

the period equations (10.9) are

$$\sum_{1\leq i,j\leq 3} P_k^{ij} x_i x_j = \overline{\sum_{1\leq i,j\leq 3} P_k^{ij} y_i y_j} \quad (k=1,3)$$

$$\sum_{1\leq i,j\leq 3} P_k^{ij} x_i y_j \in i\mathbb{R} \quad (k=1,3).$$

Now let $(i,j,k)$ be a permutation of $(1,2,3)$ and make the particular choice

$$s_1 = x_i \hat{t}_i + x_j \hat{t}_j,$$
$$s_2 = \hat{t}_k.$$

The second period equation above is satisfied for all $x_i$, $x_j$, and the first period equation can be written in the form

$$\begin{pmatrix} x_i^2 \\ x_j^2 \end{pmatrix} = \begin{pmatrix} 1 & 1 \\ e_i & e_j \end{pmatrix}^{-1} B \begin{pmatrix} 1 \\ \overline{e_k} \end{pmatrix}$$

where $\eta_i = \zeta(\omega_i)$ and $e_i = \wp(\omega_i)$ and $B$ is defined in section 27. The condition that the surface be immersed is that $s_1$ and $s_2$ have no common zeros. The zeros of $s_2$ are at $\{\omega_k/2, \omega_k/2 + \omega_1, \omega_k/2 + \omega_2, \omega_k/2 + \omega_3\}$, and

$$\hat{t}_m^2(\omega_k/2) = \hat{t}_m^2(\omega_k/2 + \omega_l) = 4(e_k - e_i) \quad (m = i, j; \; l = 1, 2, 3).$$

A necessary condition that the surface branch is that

$$(e_k - e_i)x_i^2 - (e_k - e_j)x_j^2 = 0,$$

or

$$\begin{pmatrix} g_2/2 - 3e_k^2 & -3e_k \end{pmatrix} B \begin{pmatrix} 1 \\ \overline{e_k} \end{pmatrix} = 0.$$

With the choice $\{i, j, k\} = \{1, 2, 3\}$ it can be shown that this condition is not satisfied in the standard fundamental region of the moduli space of tori. The proof uses the $q$-expansion for $g_2$ and $\eta$ given in section 27, as well as the expansion

$$e_1 = \frac{\pi^2}{6\omega_1^2}\left(1 + 24\sum_{n=1}^{\infty} \tau(n)q^n\right),$$

where

$$\tau(n) = \sum_{\substack{d|n \\ d \text{ odd}}} d.$$

Thus we have found a single immersion of every conformal type of torus punctured at the half-lattice points. Since the period conditions amount to at most six real conditions on 12 variables, there is a real 6-parameter family of surfaces, which modulo the action of the group in equation (9.8) leaves a 2-parameter family. The existence of the real two-dimensional family follows from the fact that the condition of being immersed is an open analytic condition. □

## 29 Klein bottles: conformal type, spin structure and periods

Theorem 28 shows that the torus underlying a Klein bottle must have the conformal type of the complex plane modulo a rectangular lattice, and it computes the order-two deck transformation for the covering of the Klein bottle by the torus. The theorem further shows that the torus which doubly covers the immersed Klein bottle must be untwisted. (This can also be seen from purely topological considerations.)

THEOREM 28. *Let $X : K' \longrightarrow \mathbb{R}^3$ be a complete minimal immersion of a punctured Klein bottle with finite total curvature, $\pi : T \longrightarrow K = \overline{K'}$ the oriented two-sheeted covering by a torus $T$, and $I : T \longrightarrow T$ the order-two orientation-reversing deck transformation for this cover. Then we have the following.*
  (i) *$T$ is conformally equivalent to $\mathbb{C}/\Lambda$, where $\Lambda$ is a rectangular lattice with generators $2\omega_1 \in \mathbb{R}$ and $2\omega_3 \in i\mathbb{R}$.*
 (ii) *On this torus, the deck transformation $I$ may be chosen to be $I(u) = \bar{u} + \omega_1$.*
(iii) *With this choice, the admissible spin structures are those represented by $(\wp(u) - \wp(\omega_2))du$ and $(\wp(u) - \wp(\omega_3))du$.*
(iv) *If $(s_1, s_2)$ is the spin representative of $X \circ \pi$ on $T$, the period conditions reduce to the conditions $\int_{\gamma_1} s_1^2 = 0$ and $\int_{\gamma_1} s_1 s_2 = 0$ along a closed curve $\gamma_1$ parallel to $\omega_1$.*

*Proof of (i) and (ii):* Let $\Lambda_0$ be a lattice such that $T = \mathbb{C}/\Lambda_0$. Since every conformal map from $T$ to $T$ must be linear in the standard coordinate $u$ on $\mathbb{C}$ and since $I$ is anticonformal, $I(u) = \alpha \bar{u} + \beta$ for some $\alpha, \beta \in \mathbb{C}$. The periodicity of $I$ and $I^{-1}$ implies that $\alpha \overline{\Lambda_0} \subseteq \Lambda_0$ and $\overline{\alpha}^{-1} \overline{\Lambda_0} \subseteq \Lambda_0$. These together imply that $\alpha \overline{\Lambda_0} = \Lambda_0$. Choose $\gamma \in \mathbb{C}$ satisfying $|\gamma| = 1$ and $\overline{\gamma}/\gamma = \alpha$; the rotated lattice $\Lambda = \gamma \Lambda_0$ satisfies $\overline{\Lambda} = \Lambda$ (a so-called *real* lattice). Hence $\Lambda$ is either rectangular with generators $2\omega_1 \in \mathbb{R}$, $2\omega_3 \in i\mathbb{R}$, or $\Lambda$ is rhombic with generators $2\omega_1$ and $2\omega_3 = 2\overline{\omega}_1$. On $\mathbb{C}/\Lambda$ we have $I(u) = \alpha \bar{u} + \beta$ for some new $\alpha, \beta \in \mathbb{C}$. As before, $\alpha \overline{\Lambda} = \Lambda$, but $\overline{\Lambda} = \Lambda$, so $\alpha = \pm 1$. If $\alpha = -1$, replacing $\Lambda$ by $i\Lambda$ preserves its reality, and changes $\alpha$ to 1.

52                                                                         Kusner and SchmittWith $\alpha = 1$, the condition that $I$ is involutive is that $\beta + \overline{\beta} \in \Lambda$. By the change of coordinate $u \mapsto u - i\,\text{Im}\,\beta$, it can be assumed that $\beta \in \mathbb{R}$. Then the involutive condition is that $2\beta \in \Lambda$. If $\beta \in \Lambda$ then $0$ is a fixed point of $I$. Hence $\beta \equiv \omega_1$ (rectangle) or $\beta = \omega_1 + \omega_3$ (rhombus). In the latter case, $\omega_1$ is a fixed point of $I$, so the only admissible case is the rectangle, with $I(u) = \overline{u} + \omega_1$.

*Proof of (iii):* The compatibility condition in Theorem 11 demands that $I^*I^*(s) = -s$ for any section $s$ of the spin structure. A computation shows that this condition is met only for the two spin structures named.

*Proof of (iv):* Let $\gamma_1$ and $\gamma_3$ be respectively the closed curves $t \mapsto \omega_1 t/|\omega_1| + c_1$ and $t \mapsto \omega_3 t/|\omega_3| + c_2$, $(0 \leq t \leq 2)$, where $c_1, c_2 \in \mathbb{C}$ are chosen so that the curves do not pass through any ends. Then $I(\gamma_1) = \gamma_1$, $I(\gamma_3) = -\gamma_3$. The periods conditions are

$$\int_{\gamma_k} s_1^2 = \overline{\int_{\gamma_k} s_2^2} \quad (k = 1, 3),$$

$$\int_{\gamma_k} s_1 s_2 \in i\mathbb{R} \quad (k = 1, 3).$$

With $I$ as above, under the double-cover assumption

$$(s_1, s_2) = \pm(i\overline{I^* s_2}, -i\overline{I^* s_1}),$$

we have

$$\int_{\gamma_3} s_1^2 = \int_{\gamma_3} -\overline{I^* s_2^2} = -\overline{\int_{I(\gamma_3)} s_2^2} = \overline{\int_{\gamma_1} s_2^2}$$

$$\int_{\gamma_3} s_1 s_2 = \int_{\gamma_3} \overline{I^* s_1 s_2} = \overline{\int_{I(\gamma_3)} s_1 s_2} = -\overline{\int_{\gamma_3} s_1 s_2},$$

so the period conditions are automatically satisfied for $k = 3$. Moreover, we also have

$$\int_{\gamma_1} s_1^2 = \int_{\gamma_1} -\overline{I^* s_2^2} = -\overline{\int_{I(\gamma_1)} s_2^2} = -\overline{\int_{\gamma_1} s_2^2}$$

$$\int_{\gamma_1} s_1 s_2 = \int_{\gamma_1} \overline{I^* s_1 s_2} = \overline{\int_{I(\gamma_1)} s_1 s_2} = \overline{\int_{\gamma_1} s_1 s_2}$$

and the first two period conditions (10.9) become

$$\int_{\gamma_1} s_1^2 = 0$$

$$\int_{\gamma_1} s_1 s_2 = 0$$

(this amounts to three real conditions because under the above assumption, the second integral is automatically real). □

## 30 Minimal Klein bottles with embedded planar ends

A minimal Klein bottle is constructed in this section. Its compactification is a $W$-critical surface with energy $W = 16\pi$, which lies in the *amphichiral* regular homotopy class $\mathbf{K}_0 = \mathbf{B} \# \overline{\mathbf{B}}$ of Klein bottles (cf. [13], [24]). Clearly there are no minimal Klein bottles with two embedded ends and we conjecture there are none with three embedded planar ends.

<u>Theorem</u> 29. *There exists a minimal immersion of the Klein bottle with four embedded planar ends.*

To construct this example, let $T = \mathbb{C}/\{2\omega_1, 2\omega_3\}$ be a square lattice with

$$\begin{aligned}
\omega_3 &= i\omega_1 \\
\omega_2 &= -\omega_1 - \omega_3 \\
\wp(\omega_1) &= 1 \\
\wp(\omega_2) &= 0 \\
\wp(\omega_3) &= -1.
\end{aligned}$$

Let $I: T \longrightarrow T$ be the deck transformation

$$I(u) = \bar{u} + \omega_1$$

as in Theorem 28 i). Let $a \in T$ be a point (yet to be determined) such that

$$I(a) = -a,$$

and let $E = \{a_1, \ldots, a_8\} \subset T$ be the points in Table 3. We want to construct a minimal immersion $X: (T \setminus E)/I \longrightarrow \mathbb{R}^3$,

$$X(z) = \operatorname{Re} \int^z (s_1^2 - s_2^2, i(s_1^2 + s_2^2), 2s_1 s_2),$$

where $s_1, s_2$ are sections of the spin structure $S$ determined by $\varphi$, where

$$\varphi^2 = \frac{du}{\wp(u) - \wp(\omega_2)} = \frac{du}{\wp(u)}.$$

Table 3: Values of $\wp$ and $\wp'$ at ends of Klein bottle

| $u$ | $\wp(u)$ | $\wp'(u)$ | $I(u)$ |
|---|---|---|---|
| $a_1 = a$ | $r$ | $r'$ | $a_5$ |
| $a_2 = a + \omega_2$ | $-1/r$ | $r'/r^2$ | $a_6$ |
| $a_3 = -ia$ | $-r$ | $-ir'$ | $a_4$ |
| $a_4 = -ia + \omega_2$ | $1/r$ | $-ir'/r^2$ | $a_3$ |
| $a_5 = -a_1$ | $r$ | $-r'$ | $a_1$ |
| $a_6 = -a_2$ | $-1/r$ | $-r'/r^2$ | $a_2$ |
| $a_7 = -a_3$ | $-r$ | $ir'$ | $a_8$ |
| $a_8 = -a_4$ | $1/r$ | $ir'/r^2$ | $a_7$ |

Step 1: Determination of the ends

Let $\{t_1, \ldots, t_8\}$,
$$t_\alpha = \frac{\wp(u)}{\wp(u) - \wp(a_\alpha)} \varphi$$
$$t_{\alpha+4} = \frac{\wp'(u)}{\wp(u) - \wp(a_\alpha)} \varphi \qquad (1 \le \alpha \le 4)$$

be a basis for $\mathcal{F}$, as in section 26. The skew-symmetric matrix for $\Omega$ in this basis is

$$\left(\begin{array}{c|c} 0 & W \\ \hline -W^t & 0 \end{array}\right),$$

where $W$ is given by

$$W = \begin{pmatrix} \frac{r^2+1}{r(r^2-1)} & \frac{4r}{r^2+1} & \frac{2}{r} & \frac{4r}{r^2-1} \\ \frac{-4r}{r^2+1} & \frac{r(r^2+1)}{r^2-1} & \frac{4r}{r^2-1} & -2r \\ \frac{-2}{r} & \frac{-4r}{r^2-1} & \frac{-(r^2+1)}{r(r^2-1)} & \frac{-4r}{r^2+1} \\ \frac{-4r}{r^2-1} & 2r & \frac{4r}{r^2+1} & \frac{-r(r^2+1)}{r^2-1} \end{pmatrix}.$$

The desired sections $s_1, s_2$ lie in $\ker \Omega$, so a necessary condition for existence is that

$$0 = \det W = \frac{(3r^8 - 4r^6 + 50r^4 - 4r^2 + 3)^2}{(r^4-1)^2} = \frac{9(r^4 + mr^2 + 1)^2(r^4 + \overline{m}r^2 + 1)^2}{(r^4-1)^2},$$

where $m = -2(1 - 4\sqrt{2}i)/3$. Let $r$ be the root of $r^4 + mr^2 + 1$ in the fourth quadrant; with this choice, the domain $T \backslash E$ is shown below.

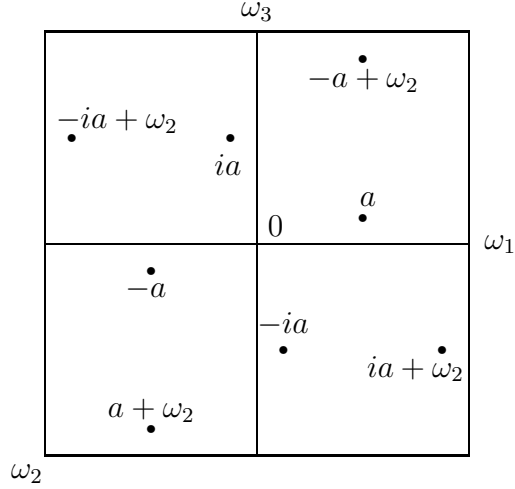

Figure 8: The eight ends in the double cover of the Klein bottle

Step 2: Choosing sections $s_1, s_2$ of $S$; the period equations

With $r$ fixed as above, rank $\Omega$ is 4, and a basis for $\ker \Omega$ is $\{\hat{s}_1, \hat{s}_2, \hat{s}_3, \hat{s}_4\}$ where

$$\hat{s}_1 = \sum_{a=1}^{4} c_1^\alpha t_\alpha$$
$$\hat{s}_2 = \sum_{a=1}^{4} c_2^\alpha t_\alpha$$
$$\hat{s}_3 = i\overline{I^* \hat{s}_1}$$
$$\hat{s}_4 = i\overline{I^* \hat{s}_2},$$

$$c_1 = (2(r^2 - 1)^2, (r^2 + 1)(r^2 - 3), (r^2 + 1)(3r^2 - 1), -2(r^2 - 1)^2)$$
$$c_2 = ((r^2 + 1)(3r^2 - 1), -2(r^2 - 1)^2, 2(r^2 - 1)^2, (r^2 + 1)(r^2 - 3))$$

and $I^*$ is a choice of a lift of the deck transformation $I$ to the spin structure $S$.

Let
$$s_1 = x_1 \hat{s}_1 + x_2 \hat{s}_2$$
$$s_2 = i\overline{I^* s_1} = \overline{x}_1 \hat{s}_3 + \overline{x}_2 \hat{s}_4.$$
$x_1, x_2 \in \mathbb{C}$

We want to find $x_1, x_2$ such that the real part of all periods are zero. By Theorem 28(iv) and section 29, the period equations reduce to the single equation

$$0 = \int_{\gamma_1} s_1^2 = x_1^2 P_1^{11} + 2x_1 x_2 P_1^{12} + x_2^2 P_1^{22},$$

where
$$P_k^{\alpha\beta} = \int_{\gamma_k} \hat{s}_\alpha \hat{s}_\beta$$

along the curve $\gamma_k : t \longmapsto t\omega_k \ (-1 \leq t \leq 1)$.

Step 3: Explicit solution of the period equation
---

The period equation above can be solved once $P_k^{\alpha\beta}$ are known. To compute these, let

$$\hat{s}_1^2 = \frac{1}{2}\left(-\sum_{\alpha=1}^{4} A_\alpha \wp(u-a_\alpha) + B\right)du, \quad A = \sum A_\alpha$$

$$\hat{s}_1 \hat{s}_2 = \frac{1}{2}\left(-\sum_{\alpha=1}^{4} C_\alpha \wp(u-a_\alpha) + D\right)du, \quad C = \sum C_\alpha$$

as in equation (B.18). Then

$$P_1^{11} = \int_{\gamma_1} \hat{s}_1^2 = A\eta_1 + B\omega_1,$$
$$P_1^{12} = \int_{\gamma_1} \hat{s}_1\hat{s}_2 = C\eta_1 + D\omega_1,$$
$$P_3^{11} = A\eta_3 + B\omega_3 = i(-A\eta_1 + B\omega_1)$$
$$P_3^{12} = C\eta_3 + D\omega_3 = i(-C\eta_1 + D\omega_1)$$
$$\eta_k = -\frac{1}{2}\int_{\gamma_k} \wp(u)du.$$

Let $J: T \to T$ be defined by $J(u) = iu$, and let $J^*$ be a lift of $J$ to $S$. Then

$$\hat{s}_1 = \sqrt{i}J^*\hat{s}_2, \quad \hat{s}_2 = \sqrt{i}J^*\hat{s}_1$$

for some choice of $\sqrt{i}$. Then

$$P_1^{12} = \int_{\gamma_1} \hat{s}_1\hat{s}_2 = \int_{\gamma_1} iJ^*\hat{s}_1\hat{s}_2 = i\int_{J_{(\gamma_1)}} J^*\hat{s}_1\hat{s}_2 = i\int_{\gamma_3} \hat{s}_1\hat{s}_2 = iP_3^{12},$$

so $D = 0$. Again,

$$P_1^{22} = \int_{\gamma_1} \hat{s}_2^2 = \int_{\gamma_1} iJ^*\hat{s}_1^2 = i\int_{J_{(\gamma_1)}} \hat{s}_1^2 = i\int_{\gamma_3} \hat{s}_1^2 = iP_3^{11},$$

so $P_1^{22} = A\eta_1 - B\omega_1$.

Having computed $P_1^{11}, P_1^{12}, P_1^{22}$ in terms of $A, B, C$, we compute $A, B, C$ by expanding $\hat{s}_\alpha \hat{s}_\beta / du$ in two ways and equating coefficients. On the one hand, by the definition of $\hat{s}_\alpha$, we have

$$\hat{s}_\alpha \hat{s}_\beta / du = \sum_{\gamma, \delta} \frac{c_\alpha^\gamma c_\beta^\delta \wp(u)}{(\wp(u) - \wp(a_\gamma))(\wp(u) - \wp(a_\delta))} \quad (1 \le \alpha, \beta \le 2; \ 1 \le \gamma, \delta \le 4).$$

Using the formula (for $\wp'(u_0)$ finite and non-zero)

$$\frac{1}{\wp(u) - \wp(u_0)} = \frac{1/\wp'(u_0)}{u - u_0} + \cdots,$$

we get the expansion at $a_\gamma$

$$\hat{s}_\alpha \hat{s}_\beta / du = \frac{c_\alpha^\gamma c_\beta^\gamma \wp(a_\gamma) / (\wp'(a_\gamma))^2}{(u - a_\gamma)^2}.$$

On the other hand we have the expansions at $a_\gamma$

$$\hat{s}_1^2 / du = \frac{-A_\gamma / 2}{(u - a_\gamma)^2}$$

$$\hat{s}_1 \hat{s}_2 / du = \frac{-C_\gamma / 2}{(u - a_\gamma)^2}.$$

Equating coefficients,

$$A_\gamma = -2(c_1^\gamma)^2 \wp(a_\gamma) / (\wp'(a_\gamma))^2$$
$$C_\gamma = -2 c_1^\gamma c_2^\gamma \wp(a_\gamma) / (\wp'(a_\gamma))^2,$$

so

$$A = \sum A_i = -32 r^2 (r^4 + 4 r^2 + 1)/3$$
$$C = \sum C_i = -2(r^4 - 1)^2.$$

To compute $B$, note that $s_1$ has a zero at 0 to get

$$B = \sum A_\gamma \wp(a_\gamma) = 4 r (r^2 + 1)^3.$$

This solves the period equation.

Finally, that the immersion is unbranched is the condition that $s_1, s_2$ have no common zeros. This amounts to the condition that if $u_0$ is a zero of $s_1$, then $I(u_0)$ is not. By using the identity

$$\overline{I^* \wp} = \frac{\wp + 1}{\wp - 1},$$

this can be checked by setting $s_1$ to zero, and solving numerically the cubic in $\wp$ which results.

# Appendix

## A  The Pfaffian

Here we recall some basic facts about skew-symmetric forms.

<u>Definition</u>.  *A bilinear form $A$ on a vector space $V$ of dimension $n$ is skew-symmetric if*

$$A(v_1, v_2) + A(v_2, v_1) = 0 \text{ for all } v_1,\ v_2 \in V,$$

*or alternatively, if the matrix $A$ for $\mathcal{A}$ satisfies*

$$A + A^t = 0.$$

The space of skew-symmetric bilinear forms is $\bigwedge^2(V^*)$. The pfaffian is a function on skew-symmetric forms whose square is the determinant.

<u>Definition</u>.  *For $A \in \bigwedge^2(V^*)$, the* pfaffian *of $A$ is*

$$pf(A) = \begin{cases} \frac{1}{m!} \overbrace{(A \wedge \ldots \wedge A)}^{m \text{ times}} & \text{if } \dim(V) = 2m \text{ is even,} \\ 0 & \text{if } \dim(V) \text{ is odd.} \end{cases}$$

For a matrix $(a_{ij})$ of $A \in \bigwedge^2(V^*)$ in the basis $\{e_1, \ldots, e_m\}$ the pfaffians for $m = 2$, $m = 4$, and $m = 6$ are respectively

$$a_{12},$$
$$a_{12}a_{34} - a_{13}a_{24} + a_{14}a_{23},$$
$$a_{12}a_{34}a_{56} - a_{12}a_{35}a_{46} + a_{12}a_{36}a_{45} - a_{13}a_{24}a_{56} + a_{13}a_{25}a_{46} -$$
$$a_{13}a_{26}a_{45} + a_{14}a_{23}a_{56} - a_{14}a_{25}a_{36} + a_{14}a_{26}a_{35} - a_{15}a_{23}a_{46} +$$
$$a_{15}a_{24}a_{36} - a_{15}a_{26}a_{34} + a_{16}a_{23}a_{45} - a_{16}a_{24}a_{35} + a_{16}a_{25}a_{34}.$$

The general pfaffian of a $2m \times 2m$ matrix has $(2m)!/(2m!) = 1 \cdot 3 \cdot 5 \cdot \ldots \cdot (2m-1)$ terms.

<u>Lemma</u>.  *The rank of a skew-symmetric matrix is even.*

*Proof.* Let $A$ be an $m \times m$ skew-symmetric matrix with rank $r$. The proof is by induction on $m$. In the case $m = 1$, then $A = (0)$ with even rank 0. Assume for some $n$ that the lemma is true for all skew-symmetric matrices smaller than $A$. If $n$ is odd, then

$$\det A = \det A^t = \det(-A) = (-1)^n \det A = -\det A,$$

so $\det A = 0$ and $A$ has a non-zero kernel. If $n$ is even, then $A$ also has a non-zero kernel unless it has full — hence even — rank $r = n$. So in either case we may assume $A$ has a non-zero kernel.

Let $v_1, \ldots, v_{n-r}$ be a basis for $\ker A$, and extend to a basis $v_1, \ldots, v_{n-r}, w_1, \ldots, w_r$ for $\mathbb{C}^n$. Let $P$ be the $n \times n$ matrix with these vectors as columns. Then $P^t A P$ is of the form

$$P^t A P = \left( \begin{array}{c|c} 0 & 0 \\ \hline 0 & A_0 \end{array} \right),$$

where $A_0$ is an $r \times r$ matrix of rank $r < n$. Moreover,

$$(P^t A P)^t = P^t A^t P = -(P^t A P),$$

so $P^t A P$, and hence $A_0$ is skew-symmetric. By the induction hypothesis, $r = \operatorname{rank} A$ is even, since it is the rank of the smaller skew-symmetric matrix $A_0$. $\square$

# B  Elliptic functions

For reference, here are some standard notations and facts about elliptic functions used in this paper (see for example [6], [7]).

*Lattices.* A non-degenerate lattice $\Lambda$ is *real* if $\Lambda = \overline{\Lambda}$. There are two kinds of real lattices:
 (i) rectangular: generators $\omega_1 \in \mathbb{R}$ and $\omega_3 \in i\mathbb{R}$ can be chosen for $\Lambda$.
 (ii) rhombic: generators $\omega_1$ and $\omega_3 = \overline{\omega}_1$ can be chosen for $\Lambda$.
For any lattice with generators $\omega_1, \omega_3$, let $\omega_2 = -\omega_1 - \omega_3$.

*The Weierstrass $\wp$ function:* Given a lattice $\Lambda$ generated by $\omega_1$ and $\omega_3$, the elliptic function $\wp$ on $\mathbb{C}/\Lambda$ satisfies the differential equation

$$(\wp')^2 = 4\wp^3 - g_2\wp - g_3 = 4(\wp - e_1)(\wp - e_2)(\wp - e_3),$$

where

$$\begin{aligned} e_i &= \wp(\omega_i) \quad (i = 1, 2, 3), \\ e_1 + e_2 + e_3 &= 0, \\ g_2 &= -4(e_1 e_2 + e_1 e_3 + e_2 e_3), \\ g_3 &= 4 e_1 e_2 e_3. \end{aligned}$$

The function $\wp$ has a double pole at $0$ and two simple zeros which come together only on the square lattice; $\wp'$ has a triple pole at $0$ and three simple poles at $\omega_1, \omega_2, \omega_3$.

The function $\wp$ is even; $\wp'$ is odd. On a horizontal rectangular lattice, $\wp(\overline{u}) = \overline{\wp(u)}$; on a horizontal square lattice, $\wp(iu) = -\wp(u)$.

The expansion for $\wp$ at 0 is

$$\wp(u) = \frac{1}{u^2} + \frac{g_2}{20}u^2 + \dots.$$

A useful property of $\wp$ is the following special case of the addition formula ($\{i, j, k\}$ is any permutation of $\{1, 2, 3\}$):

(B.15) $$\wp(u \pm \omega_i) = e_i + \frac{(e_i - e_j)(e_i - e_k)}{\wp(u) - e_i}.$$

*The Weierstrass $\zeta$ function:* The $\zeta$ function is defined by

$$\zeta(u) = -\int \wp(u) du,$$

with the constant of integration chosen so that $\lim_{u \to 0} \zeta(u) - u^{-1} = 0$. With $\eta_i = \zeta(\omega_i)$ ($i = 1, 2, 3$), properties of $\zeta$ include:

$$\eta_1 + \eta_2 + \eta_3 = 0,$$
$$\zeta(u + 2\omega_i) = \zeta(u) + 2\eta_i \quad (i = 1, 2, 3),$$
$$\zeta \text{ is an odd function.}$$

Legendre's relation is that

(B.16) $$\eta_1 \omega_3 - \eta_3 \omega_1 = i\pi/2.$$

A form of the quasi-addition formula for $\zeta$ is

(B.17) $$\zeta(u - v) - \zeta(u) + \zeta(v) = \frac{1}{2}\left(\frac{\wp'(u) + \wp'(v)}{\wp(u) - \wp(v)}\right).$$

A useful property of elliptic functions which can also be stated in more generality is the following: Let $f$ be an elliptic function with poles of order at most 2, with no residues, and with principal parts

$$\frac{a_1}{(u - \alpha_1)^2}, \dots, \frac{a_n}{(u - \alpha_n)^2}.$$

Then

(B.18) $$f(u) = b + \sum a_i \wp(u - a_i)$$

for some $b$, because the difference $f(u) - \sum \alpha_i \wp(u - \alpha_i)$ has no poles and hence is constant.

# References

[1] Arbarello, E., Cornalba, M., Griffiths, P. A., and Harris, J. *Geometry of algebraic curves.* New York: Springer-Verlag, 1985.

[2] Bryant, R. *A duality theorem for Willmore surfaces.* J. Diff. Geom., 20:23–53, 1984.

[3] Bryant, R. *Surfaces in conformal geometry.* Proceedings of Symposia in Pure Mathematics 48:227–240, 1988.

[4] Callahan, M., Hoffman, D., Meeks III, W. H. *Embedded minimal surfaces with an infinite number of ends.* Invent. Math. 96:459–505, 1989.

[5] Costa, C. *Complete minimal surfaces in $\mathbb{R}^3$ of genus one and four planar embedded ends.* Preprint, 1990.

[6] DuVal, P. *Elliptic functions and elliptic curves.* Cambridge: Cambridge University Press, 1973.

[7] Erdelyi, A., ed. *Higher transcendental functions.* New York: McGraw-Hill, 1953.

[8] Gilbert, J. E., and Murray, A. M. A. *Clifford algebras and Dirac operators in harmonic analysis.* Cambridge: Cambridge University Press, 1991.

[9] Griffiths, P., and Harris, J. *Principles of algebraic geometry.* New York: Wiley-Interscience, 1978.

[10] Gunning, R. *Lectures on Riemann surfaces.* Princeton: Princeton University Press, 1966.

[11] Gunning, R. *Riemann surfaces and generalized theta functions.* New York: Springer Verlag, 1976.

[12] Kauffman, L. H. *On knots.* Princeton: Princeton University Press, 1987.

[13] Kusner, R. *Comparison surfaces for the Willmore problem.* Pacific Journal of Mathematics 138:317–345, 1989.

[14] Kusner, R. *Conformal geometry and complete minimal surfaces.* Bull. Amer. Math. Soc. 17:291–295, 1987.

[15] Kusner, R. *Global geometry of extremal surfaces in three-space.* Dissertation, Univ. of California, Berkley, 1988.

[16] Lang, S. *Elliptic functions.* New York: Springer-Verlag, 1987.

[17] Lawson, H. B., and Michaelsohn, M. L. *Spin Geometry* Princeton: Princeton U. Press, 1989.

[18] Meeks III, W. H. *The classification of complete minimal surfaces in $\mathbb{R}^3$ with total curvature greater than $-8\pi$.* Duke Mathematical Journal 48:523–535, 1981.

[19] Meeks III, W. H. and Patrusky, J. *Representing homology classes by embedded circles on a compact surface.* Ill. J. Math. 22:262-269, 1978.

[20] Milnor, J. *Spin structures on manifolds.* Enseign. Math. 9:198–203, 1963.

[21] Mumford, D. *Tata lectures on theta,* v. 1 and 2. Boston: Birkhäuser, 1983.

[22] Osserman, R. *A survey of minimal surfaces.* New York: Van Nostrand Reinhold, 1969.

[23] Peng, C. K. *Some new examples of minimal surfaces in $\mathbb{R}^3$ and its applications.* Preprint, MSRI, 1986.

[24] Pinkall, U. *Regular homotopy classes of immersed surfaces.* Topology, 24:421–434, 1985.

[25] Rees, E. *Notes on geometry.* New York: Springer-Verlag, 1983.

[26] Schmitt, N. *Minimal surfaces with embedded planar ends.* Dissertation, Univ. of Massachusetts, Amherst, 1993.

[27] Sullivan, D. *The spinor representation of minimal surfaces in space.* Notes, 1989.

# Graphics

The following computer graphics were produced at the Center for Geometry, Analysis, Numerics and Graphics using the MESH program authored by Jim Hoffman. For more information about MESH, contact `web@gang.umass.edu`.